\documentclass[12pt]{amsart}

\usepackage{amsmath}
\usepackage{amsthm}
\usepackage{amssymb}
\usepackage{graphicx}

\begingroup 

\newtheorem{theorem}{Theorem}[section]  
\newtheorem{corollary}[theorem]{Corollary}     
\newtheorem{lemma}[theorem]{Lemma}         
\newtheorem{proposition}[theorem]{Proposition}  
\endgroup

\newtheorem{definition}[theorem]{Definition}   

\newtheorem{remark}[theorem]{Remark}        
\newtheorem{example}[theorem]{Example}        

\numberwithin{equation}{section}

\newcommand{\R}{{\mathbb R}}
\newcommand{\N}{{\mathbb N}}

\newcommand{\e}{\varepsilon}

\newcommand{\vp}{\varphi}

\newcommand{\diam}{\operatorname{diam}}

\newcommand{\mcA}{{\mathcal{A}}}
\newcommand{\mcB}{{\mathcal{B}}}

\newcommand{\mcJ}{{\mathcal{J}}}

\newcommand{\mcM}{{\mathcal{M}}}
\newcommand{\mcP}{{\mathcal{P}}}

\begin{document}

\title[Generalized $\alpha$-variation and differentiable functions]{Generalized $\alpha$-variation and Lebesgue equivalence to
differentiable functions}

\author{Jakub Duda}
\thanks{The author was supported in part by ISF}
\address{
Department of Mathematics, Weizmann Institute of Science,
Rehovot 76100, Israel}
\curraddr{Seifertova 19, 130 00 Prague 3, Czech Republic}
\email{jakub.duda@gmail.com}

\date{\today}


\begin{abstract} 
We find an equivalent condition for a real function $f:[a,b]\to\R$ to be Lebesgue equivalent to
an $n$-times differentiable function ($n\geq 2$); a simple solution in the case $n=2$ appeared
in an earlier paper.
For that purpose, we introduce the notions of  $CBVG_{1/n}$ and $SBVG_{1/n}$ functions, which play analogous r\^oles
for the $n$-th order differentiability as the classical notion of a $VBG_*$ function
for the first order differentiability, and the classes $CBV_{1/n}$ and $SBV_{{1}/{n}}$
(introduced by Preiss and Laczkovich) for $C^n$ smoothness.
As a consequence of our approach, we obtain that Lebesgue equivalence to $n$-times
differentiable function is the same as Lebesgue equivalence to a function $f$ which
is $(n-1)$-times differentiable with $f^{(n-1)}(\cdot)$ being pointwise Lipschitz.
We also characterize the situation when a given function is Lebesgue equivalent
to an $n$-times differentiable function $g$ such that $g'$ is nonzero a.e.
As a corollary, we establish a generalization of Zahorski's Lemma for higher
order differentiability.
\end{abstract}

\maketitle

\section{Introduction} 

Let $f:[a,b]\to\R$. We say that {\em $f$ is Lebesgue equivalent to $g:[a,b]\to\R$}
provided there exists a homeomorphism $h$ of $[a,b]$ onto itself such that $g=f\circ h$.
This terminology is taken from~\cite{Cesari}.
Zahorski~\cite{Zah1} and Choquet~\cite{Ch} (see also Tolstov~\cite{Tol1})
proved a result
characterizing paths ($f:[a,b]\to\R^n$) that
allow a differentiable parametrization (resp.\ a dif.\ parametrization with almost everywhere nonzero
derivative) as those paths that have
the $VBG_*$ property (resp.\ which are also not constant on any interval). 
Fleissner and Foran~\cite{FF} reproved this later (for real functions only and not considering
the case of a.e.\ nonzero derivatives) using a different result
of Tolstov.
The definition of $VBG_*$ is classical; see e.g.~\cite{S}.
The mentioned results were generalized by L.~Zaj\'\i\v{c}ek and the author~\cite{DZ}
to paths with values in Banach spaces 
(and also metric spaces using the metric derivative instead
of the usual one).
Laczkovich, Preiss~\cite{LP}, and Lebedev~\cite{Leb}
studied (among other things) the case of $C^n$-parametrizations of
real-valued functions ($n\geq2$).
Lebedev proved that a continuous function $f:[a,b]\to\R$ is Lebesgue
equivalent to a $C^n$ function provided
\[\lambda(f(K_f))=0\text{ and }\sum_{\alpha\in A}(\omega^f_\alpha)^{1/n}<\infty,\]
where $K_f$ is the set of point of varying monotonicity of~$f$
(see the definition below) and $\omega^f_\alpha$ is the oscillation of $f$
on $I_\alpha$, where $(I_\alpha)_{\alpha\in A}$ are all the intervals contiguous
to~$K_f$ in $[a,b]$.
Laczkovich and Preiss showed that the same is true for a continuous $f$
provided 
\begin{equation}\label{lpv}
V_{1/n}(f,K_f)<\infty,
\end{equation}
or 
\begin{equation}\label{lpsv}
SV_{1/n}(f,K_f)=0.
\end{equation}
(See Definition~\ref{lpdef} in Section~\ref{prelim}).
They define $CBV_{1/n}$ (resp.\ $SBV_{1/n}$)
as the class of continuous function which satisfy~\eqref{lpv}
(resp.\ \eqref{lpsv}).
Moreover, in~\cite{LP} and~\cite{Leb} also the case of $C^{n,\alpha}$ ($0<\alpha\leq1$)
parametrizations is settled (where $C^{n,\alpha}$ is the class
of functions such that $f^{(n)}$ is $\alpha$-H\"older).
\par
Differentiability via a homeomorphic change of variable
was studied by other authors (see e.g.~\cite{BG},~\cite{Br}).
For a nice survey of differentiability of real-valued
functions via homeomorphisms, see~\cite{GNW}.
L.~Zaj\'\i\v{c}ek and the author~\cite{DZC2,DZDC} characterized
the situation when a Banach space-valued path (for Banach spaces with a $C^1$ norm) 
admits a $C^2$-parametrization or a parametrization with finite convexity.
In the corresponding situations, also the case of the first derivative
being almost everywhere nonzero is treated in~\cite{DZC2,DZDC}.
\par
In~\cite{D12}, we characterized the situation when a function $f:[a,b]\to\R$ 
is Lebesgue equivalent to a twice differentiable function.
We introduced the notion of $VBG_{1/2}$ functions for that purpose. 
We also established that for a real function $f$ defined on a closed interval,
being Lebesgue equivalent to a twice differentiable
function is equivalent to being Lebesgue equivalent to a differentiable function
whose derivative is pointwise Lipschitz.
\par
In the present article, we characterize the situation when $f$ is Lebesgue
equivalent to an $n$-times differentiable function for $n\geq3$
(our approach in the present article gives a certain condition also in case $n=2$ which
can be seen to be equivalent to the one proved in~\cite{D12},
but the present general proof is much more complicated than the arguments of~\cite{D12} 
in that interesting special case).
We introduce two
new classes of functions: $CBVG_{1/n}$ and $SBVG_{1/n}$, which are analogous 
to the classes $CBV_{1/n}$ and $SBV_{1/n}$ introduced by Preiss and Laczkovich
in~\cite{LP} in order to characterize the situation when a given function is Lebesgue
equivalent to a $C^n$ function. 
In the main Theorem~\ref{mainthm}, we prove
that $f$ is Lebesgue equivalent to an $n$-times differentiable function
if and only if $f$ is $CBVG_{1/n}$ (resp.\ $f$ is $SBVG_{1/n}$).
As a corollary, we obtain that the classes $CBVG_{1/n}$ and $SBVG_{1/n}$ coincide
(which seems to be difficult to establish directly).
Our approach also
yields that $f$ is Lebesgue equivalent to an $n$-times differentiable function  
if and only if $f$ is Lebesgue equivalent
to an $(n-1)$-times differentiable function $g$ such that $g^{(n-1)}(\cdot)$ is pointwise Lipschitz
(see Theorem~\ref{mainthm}).
This corresponds to the analogous situation for $n=2$ in~\cite{D12}, and is a similar
phenomenon as we see with $f:[a,b]\to\R$ being Lebesgue equivalent to a $C^n$ function
if and only if $f$ is Lebesgue equivalent to a $C^{(n-1),1}$ function or
a function which has bounded $n^{\text{th}}$ derivative (see~\cite[Remark~3.7]{LP}), which is proved
in~\cite{LP}.
We also present an example (see Example~\ref{difex}) which shows that
for each $n\geq2$ there exists a continuous function, which is
$CBVG_{1/n}$, but which is not Lebesgue equivalent to any $C^n$ function.
In Theorem~\ref{mainnonzerothm}, we characterize Lebesgue equivalence
to an $n$-times differentiable function whose first derivative is a.e.\ nonzero.
\par
The classical Zahorski's Lemma (see e.g.~\cite{Zah1}
or~\cite[p.~27]{GNW}) claims that if $M\subset[a,b]$ has Lebesgue measure~$0$,
then there exists a (boundedly) differentiable homeomorphism~$h$ from $[a,b]$ onto itself such
that $h^{-1}(M)\subset\{x\in[a,b]:h'(x)=0\}$. If $M$ is closed,
then $h$ can be taken~$C^1$. In Theorem~\ref{zahthm}, we show a higher
order analogue of this fact; i.e., a closed set $M\subset[a,b]$
is the image by an $n$-times differentiable homeomorphism such that $h^{(i)}(x)=0$
for all $x\in h^{-1}(M)$ and $i=1,\dots,n$, if and only if there exists
a decomposition of the set $M$ such that certain variational conditions
closely related to the definition of the class
$CBVG_{1/n}$ (respectively, $SBVG_{1/n}$) are satisfied. See Theorem~\ref{zahthm} for details.
\par
The current paper is structured as follows. Section~\ref{prelim} contains  basic facts and
definitions. Section~\ref{gvarsec} contains facts about the generalized variation $GV_{1/n}$ 
(and related notions) and
classes $CBVG_{1/n}$, resp.\ $SBVG_{1/n}$; there is also a definition of an auxiliary class
$\overline{SBVG}_{1/n}$. 
Section~\ref{mainsec} contains the main Theorems~\ref{mainthm} and~\ref{mainnonzerothm}.
In section~\ref{Zahorski} we prove Theorem~\ref{zahthm}, which
is an analogue of the Zahorski's Lemma 
for higher order differentiability.
\par
In the proofs, we need many auxiliary results. Let us point out that the main 
ingredients for our results are the estimate of Lemma~\ref{nn1odhlem},
and the method of construction of a suitable variation in Lemmata~\ref{propvlem},~\ref{proptildevlem}.

\section{Preliminaries}\label{prelim}

By $C$ (resp.\ $C_x$, \dots) we will denote an absolute constant (resp.\ constant depending
on~$x$, \dots) that 
can change between lines. By letter $n$ we will always denote a positive integer.
By $\lambda$ we will denote the Lebesgue measure on $\R$.
For $x,r\in\R$ with $r>0$ we will denote by $B(x,r):=\{y\in\R:|x-y|<r\}$
the open ball with center $x$ and radius~$r$.
Let $K\subset [a,b]$ be closed, and such that $\{a,b\}\subset K$. 
We say that the interval $(c,d)\subset[a,b]$ is {\em contiguous to $K$ $($in $[a,b])$}
provided $c,d\in K$ and $(c,d)\cap K=\emptyset$ (i.e.\ it is a maximal open
component of $[a,b]\setminus K$ in $[a,b]$).
Let $f:[a,b]\to\R$. 
By $K_f$ we will denote {\em the set of points of varying monotonicity of $f$},
i.e.\ the set of points $x\in[a,b]$ such that there is no open neighbourhood $U$ of $x$
such that $f|_{U}$ is either constant or strictly monotone (see e.g.\ \cite{LP}).
Obviously, $K_f$ is closed and $\{a,b\}\subset K_f$.
We will also frequently use the simple fact that if $h$ is a homeomorphism of $[a,b]$ onto itself, $g=f\circ h$,
then $K_g=h^{-1}(K_f)$.
\par
Let $f:[a,b]\to\R$. We say that $f$ is {\em pointwise Lipschitz at $x\in[a,b]$}
provided 
\[ \limsup_{\substack{t\to0\\x+t\in[a,b]}} \frac{|f(x+t)-f(x)|}{|t|}<\infty.\]
We say that $f$ is {\em pointwise Lipschitz} provided it is pointwise Lipschitz
at each point $x\in[a,b]$.
We will define the derivative $f'(x)$ of $f$ at $x\in[a,b]$ as usual;
at the endpoints we consider the corresponding unilateral derivative.
The $n^{\text{th}}$ derivative $f^{(n)}(x)$ of $f$ at $x$ is defined by 
induction as $f^{(0)}(x):=f(x)$, and $f^{(n)}(x):=(f^{(n-1)})'(x)$ for $n\geq1$.
We say that $f$ is {\em $C^n$} for $n\geq1$ provided $f^{(n)}$ exists and is continuous in $[a,b]$.
We will often use the following easy fact: if $f$ is $C^1$, and $x\in K_f\cap(a,b)$, then $f'(x)=0$.

The following version of Sard's theorem is proved in e.g.\ \cite[Lemma~2.2]{DZ}.

\begin{lemma}\label{Sard} 
If $f:\R\to\R$, then $\lambda(f(\{x\in\R:f'(x)=0\}))=0$.
\end{lemma}

The following simple lemma is proved in~\cite[Lemma~9]{D12}.

\begin{lemma}\label{propnlem}
Let $h_m:[a,b]\to[c_m,d_m]$ $(m\in\mcM\subset\N)$ be continuous increasing
functions such that $\sum_{m\in\mcM}h_m(x)\in\R$ for all $x\in[a,b]$.
Let $K\subset[a,b]$ be closed and such that $\lambda(h_m(K))=0$ for all $m\in\mcM$.
Then $h:[a,b]\to[c,d]$, defined as $h(x):=\sum_{m\in\mcM}h_m(x)$, 
is a continuous and increasing function $($for some
$c,d\in\R)$ such that $\lambda(h(K))=0$.
\end{lemma}

The following definition is taken from~\cite{LP}.

\begin{definition}\label{lpdef}
Let $f:[a,b]\to\R$ be continuous. For $\alpha\in(0,1]$, $\delta>0$, and
$K\subset[a,b]$, we will define $V_{\alpha}^\delta(f,K)$ as a supremum of the
sums
\begin{equation}\label{v1overm} 
\sum^N_{i=1}|f(d_i)-f(c_i)|^{\alpha},
\end{equation}
where $([c_i,d_i])_{i=1}^N$ is any sequence of
non-overlapping intervals in $[a,b]$ such that $c_i,d_i\in K$, $d_i-c_i\leq\delta$ for
all $i=1,\dots,N$. We define $V(f,[a,b]):=V_{1}^{b-a}(f,[a,b])$,
$V_{\alpha}(f,K):=V_{\alpha}^{b-a}(f,K)$,
and 
\[SV_{\alpha}(f,K):=\lim_{\delta\to0+}V^\delta_{\alpha}(f,K).\]
\end{definition}

See the paper~\cite{LP} for basic properties of these fractional variations.
Note that~\cite[Lemma~3.13]{LP} implies
that if $K\subset[a,b]$ is closed, and $V_{\alpha}^\delta(f,K)<\infty$, then 
the function $g(x)=V_{\alpha}^\delta(f,K\cap[a,x])$ is continuous on~$K$.

\par
We have the following simple lemma.

\begin{lemma}\label{zerolem}
Suppose that $H\subset\R$ is bounded, $f:H\to\R$ is uniformly continuous
and such that $V_{\alpha}^\delta(f,H)<\infty$ for some $\alpha\in(0,1)$ and $\delta>0$.
Then $\lambda(\overline{f(H)})=0$.
\end{lemma}

\begin{proof} 
Let $\e>0$ and put $\eta:=\big(\frac{\e}{V_\alpha^\delta(f,H)+1}\big)^{\frac{1}{1-\alpha}}$. 
Choose $0<\zeta<\delta$ such that for
all $x,y\in H$ with $|x-y|<\zeta$ we have $|f(x)-f(y)|<\eta$.
Let $([a_i,b_i])_{i=1}^N$ be non-overlapping intervals with
$a_i,b_i\in H$ and $b_i-a_i<\zeta$. Then 
\[\sum_{i=1}^N |f(b_i)-f(a_i)|\leq \eta^{(1-\alpha)}\cdot V^\delta_{\alpha}(f,H)<\e.\]
Thus we have $SV_1(f,H)=0$;  
\cite[Theorem~2.9]{LP} implies $\lambda(\overline{f(H)})=0$.
\end{proof}

We shall need the following lemma. For a proof, see e.g.\ \cite[Lemma~2.7]{DZ}.

\begin{lemma}\label{basiclem} 
Let $\{a,b\}\subset B\subset[a,b]$ be closed, and $f:[a,b]\to\R$ be
continuous. If $\lambda(f(B))=0$,
then $V(f,[a,b])=\sum_{i\in{\mathcal I}}V(f,[c_i,d_i])$,
where $I_i=(c_i,d_i)$, $(i\in{\mathcal I}\subset\N)$
are all intervals contiguous to $B$ in $[a,b]$.
\end{lemma}

\begin{lemma}\label{variacelem}
Let $f:[a,b]\to\R$ be continuous, $\{a,b\}\subset K\subset[a,b]$ be closed, 
such that $V_{\alpha}^\delta(f,K)<\infty$
for some $\alpha\in(0,1)$, $\delta>0$, and $V(f,[c,d])=|f(d)-f(c)|$
whenever $(c,d)$ is an interval contiguous to $K$ in~$[a,b]$. 
Then $V(f,[a,b])<\infty$.
\end{lemma}

\begin{proof}
Let $(u_p,v_p)$ ($p\in\mcP\subset\N$) 
be all the intervals contiguous to $K$ in $[a,b]$.
By Lemma~\ref{zerolem}, we have that $\lambda(f(K))=0$, and thus
by Lemma~\ref{basiclem} and the assumptions, we obtain
\[ V(f,[a,b])=\sum_{p\in\mcP}V(f,[u_p,v_p])
=\sum_{p\in\mcP}|f(v_p)-f(u_p)|\leq 2NM+V_{\alpha}^\delta(f,K)<\infty,\]
where $N$ is the number of $p\in\mcP$ such
that either $|f(v_p)-f(u_p)|>1$ or $v_p-u_p\geq\delta$ 
(which is finite since $V^\delta_\alpha(f,K)<\infty$ and $b-a<\infty$),
and $M:=\max_{x\in[a,b]}|f(x)|$.
\end{proof}

\begin{lemma}\label{v1nnul}
Let $f:[a,b]\to\R$ be continuous, $\{a,b\}\subset K\subset[a,b]$ closed,
$V_{\alpha}^\delta(f,K)<\infty$
for some $\alpha\in(0,1)$,  $\delta>0$, and $V(f,[u,v])=|f(v)-f(u)|$ whenever
$(u,v)$ is an interval contiguous to~$K$.
Let $g(x):=V_{\alpha}^\delta(f,K\cap[a,x])$.
Then $\lambda(g(K))=0$.
\end{lemma}

\begin{proof}
Lemma~\ref{zerolem} shows that $\lambda(f(K))=0$, and by Lemma~\ref{variacelem} it follows that $V(f,[a,b])<\infty$.
Let $\tilde g(x)=g(x)$ for $x\in K$. 
Now,~\cite[Lemma~3.13]{LP} shows that $g$ is continuous on~$K$. 
Extend $\tilde g$ to $[a,b]$
in such a way that $\tilde g$ is affine and continuous on every $[u,v]$,
whenever $(u,v)$ is an interval contiguous to $K$.
Since $g(K)=\tilde g(K)$, it is enough to prove that $\lambda(\tilde g(K))=0$.
Let $\e>0$ and $([c_i,d_i])^N_{i=1}$ be non-overlapping intervals with $c_i,d_i\in K$, $d_i-c_i\leq\delta$ for each $i=1,\dots,N$, and
such that $V^\delta_{\alpha}(f,K)\leq\sum^N_{i=1}|f(d_i)-f(c_i)|^{\alpha}+\e/2$.
For each $i=1,\dots,N$, by Lemma~\ref{basiclem} find intervals $([c^i_j,d^i_j])^{J_i}_{j=1}$ 
contiguous to $K$ in $[c_i,d_i]$ such that 
\[ |f(d_i)-f(c_i)|\leq V(f,[c_i,d_i])
\leq\sum_{j=1}^{J_i} |f(d^i_j)-f(c^i_j)|+\bigg(\frac{\e}{2N}\bigg)^{1/\alpha}.\]
Thus
\[ \tilde g(b)-\tilde g(a)=V_{\alpha}^\delta(f,K)\leq\sum_{i=1}^N\sum_{j=1}^{J_i}
|f(d^i_j)-f(c^i_j)|^{\alpha}+\e
\leq\sum_{p\in\mcP}|f(v_p)-f(u_p)|^{\alpha}+\e,\]
where $(u_p,v_p)$ ($p\in\mcP\subset\N$) are all the intervals contiguous
to $K$ in $[a,b]$.
Now, send $\e\to0+$ to conclude that 
$\tilde g(d)-\tilde g(c)\leq\sum_{p\in\mcP}\tilde g(v_p)-\tilde g(u_p)$.
Since $\tilde g(K)\cap \tilde g\big(\bigcup_{p\in\mcP}(u_p,v_p)\big)$ is countable,
we obtain that 
\[\lambda(\tilde g(K))=(\tilde g(b)-\tilde g(a))-\sum_{p\in\mcP}(\tilde g(v_p)-\tilde g(u_p))=0.\]
\end{proof}

\begin{remark}
Lemmata~\ref{propnlem} and~\ref{v1nnul} together with the methods from~\cite{LP} can be used to give a proof of the following:
{\em Suppose that $f:[a,b]\to\R$ is continuous, $n\geq2$. Then the following are equivalent:
\begin{enumerate}
\item $f$ is Lebesgue equivalent to a $C^n$ function $g$ such that $g'(x)\neq0$ for almost all $x\in[a,b]$;
\item $V_{1/n}(f,K_f)<+\infty$ and $f$ is not constant on any interval;
\item $SV_{1/n}(f,K_f)=0$ and $f$ is not constant on any interval.
\end{enumerate}
}
\noindent
See Definition~\ref{lpdef} for the definitions of the fractional variations.
In~(i), we can actually replace ``$C^n$'' with ``$C^s$'', where $s>1$ 
$($we did not define this class$)$, if we replace $1/n$ by $1/s$ in the variations in~(ii) and~(iii).
See~\cite[Definition~3.1]{LP} for the definition of this class.
\end{remark}

Let $K\subset\R$ be closed. As usual, by $K'$, we will denote the set of 
all accumulation points of~$K$.
It is easy to see that $K'$ is closed, and $K\setminus K'$ is countable. 
\par
We have the following easy consequence of the classical Rolle's theorem.

\begin{lemma}\label{Rolle}
Let $f:[a,b]\to\R$ be $(n-1)$-times differentiable for some $n\geq2$,
and suppose that there are $n$ distinct points $(x_i)^n_{i=1}$ in $[a,b]$ such
that $f'(x_i)=0$ for $i=1,\dots, n$. Then there exists $y\in[\min_{1\leq i\leq n} x_i,\max_{1\leq i\leq n} x_i]$
such that $f^{(n-1)}(y)=0$.
\end{lemma}

The following lemma shows that the derivatives are zero at all accumulation points of a given set.

\begin{lemma}\label{Kflem}
Let $f:[a,b]\to\R$ be $(n-1)$-times differentiable for some $n\geq2$, let $K\subset\{x\in[a,b]:f'(x)=0\}$
be a closed set. 
Then $f^{(i)}(x)=0$ whenever $x\in K'\cap(a,b)$ and $i\in\{1,\dots,n-1\}$.
If $x\in K'\cap(a,b)$, and $f^{(n)}(x)$ exists, then $f^{(n)}(x)=0$.
\end{lemma}

\begin{proof} By assumptions, we have that $f'(x)=0$ for all $x\in K$.
Let $i\in\{2,\dots,n-1\}$.
Without any loss of generality, assume that $x$ is a right-hand-side
accumulation point of $K$.
Fix $k\in\N$. There are $x<z_1<\dots<z_{i}<x+1/k$ with $z_j\in K$.
By~Lemma~\ref{Rolle} there exists $w_k\in[z_1,z_{i}]$ with $f^{(i-1)}(w_k)=0$,
and thus (by induction) we have 
$f^{(i)}(x)=\lim_{k\to\infty} \frac{f^{(i-1)}(w_k)}{w_k-x}=0$
for $i=2,\dots,n-1$.
\par
If $f^{(n)}(x)$ exists, then the argument above implies that it is equal to $0$.
\end{proof}

Next lemma will allow us to construct suitable extensions
of functions.

\begin{lemma}\label{homeolem}
Let $\alpha,\beta,A,B\in\R$, with $\alpha<\beta$, $A<B$, $n\in\N$, $n\geq2$.
Then there exists an $n$-times differentiable 
homeomorphism $H:[\alpha,\beta]\to[A,B]$ such that
\begin{enumerate}
\item $H^{(i)}(\alpha)=H^{(i)}(\beta)=0$ for $i=1,\dots,n$, $H'(x)>0$ for all $x\in(\alpha,\beta)$,
\item $H(\alpha)=A$, $H(\beta)=B$,
\item $|H^{(n-1)}(t)|\leq C\cdot\frac{|B-A|}{(\beta-\alpha)^n}\cdot \min(t-\alpha,\beta-t)$ for $t\in(\alpha,\beta)$,
where $C>0$ is an absolute constant.
\end{enumerate}
\end{lemma}

\begin{proof}
As in~\cite[p.\ 420]{LP}, denote
\[ w(x)=c\int^x_0\exp(-t^{-2}-(1-t)^{-2})\,dt,\quad (x\in[0,1]),\]
where the positive constant $c$ is chosen such that $w(1)=1$.
Then $w\in C^\infty([0,1])$, $w(0)=0$, $w$ is strictly increasing on~$[0,1]$, $w'(x)\neq0$
for all $x\in(0,1)$, and $w^{(i)}(0)=w^{(i)}(1)=0$ for $i=1,2,\dots$.
For $x\in[\alpha,\beta]$  define
\[ H(x)=A+(B-A)\cdot w\bigg(\frac{x-\alpha}{\beta-\alpha}\bigg).\]
It is easy to see that the conditions~(i) and~(ii) hold. Condition~(iii) follows
from the fact that the function~$w$ is $C^\infty$ and thus $H^{(n-1)}(\cdot)$ is $C\cdot\frac{|B-A|}{(\beta-\alpha)^n}$-Lipschitz
on~$[\alpha,\beta]$.
\end{proof}

\section{Generalized fractional variation}\label{gvarsec}

We need the following generalized fractional variation.

\begin{definition}\label{gvndef}
Let $f:[a,b]\to\R$. Let $\emptyset\neq A\subset K$ be closed sets.
We define the generalized $(1/n)$-variation $GV_{1/n}(f,A,K)$ 
$($resp.\ $\overline{GV}_{1/n}^\delta(f,A,K)$ for $\delta>0)$ as
the supremum of sums
\begin{equation}\label{Wn} 
\sum^{N}_{i=1} \left(V_{\frac{1}{n-1}}(f,K\cap[x_i,y_i])\right)^{\frac{n-1}{n}},
\end{equation}
where the supremum is taken over all
collections of non-overlapping intervals $([x_i,y_i])_{i=1}^N$ in $[a,b]$
with $x_i,y_i\in A$ for all $i=1,\dots,N$ $($resp.\ over all collections
$([x_i,y_i])^N_{i=1}$ of
non-overlapping intervals in $[a,b]$ such that $y_i-x_i\leq\delta$, $x_i,y_i\in K$,
and $\{x_i,y_i\}\cap A\neq\emptyset$ for all $i=1,\dots,N)$.
\par 
We put $GV_{1/n}(f,\emptyset,K)=\overline{GV}^\delta_{1/n}(f,\emptyset,K)=0$.
Similarly, we also define auxiliary variation $rGV^\delta_{1/n}(f,A,K)$ 
$($resp.\ $lGV^\delta_{1/n}(f,A,K))$ as a supremum
of the sums in~\eqref{Wn} for $\overline{GV}^\delta_{1/n}$, 
but requiring that $x_i\in A$ $($resp.\ $y_i\in A)$ for all $i=1,\dots,N$,
whenever $([x_i,y_i])_{i=1}^N$ is a~sequence of admissible intervals 
for $\overline{GV}^\delta_{1/n}$ in~\eqref{Wn}.
\par
In all cases, when there is no admissible sequence $([x_i,y_i])_{i=1}^N$,
we define the corresponding variation to be equal to~$0$.
\end{definition}

If $A\subset K\subset[a,b]$ are closed sets, $rGV^\delta_{1/n}(f,A,K)<\infty$ (respectively, $lGV^\delta_{1/n}(f,A,K)<\infty$),
and $x\in A$, then $rGV^\delta_{1/n}$ (resp.\ $lGV^\delta_{1/n}$) is ``additive at $x$'', i.e.
\begin{multline}\label{rAdd} rGV^\delta_{1/n}(f,A,K)\\=rGV^\delta_n(f,A\cap[a,x],K\cap[a,x])+rGV^\delta_n(f,A\cap[x,b],K\cap[x,b]);
\end{multline}
(and similarly for $lGV^\delta_n$).
This is easily seen from the definition. 
Also, we have that
\begin{multline}\label{basicpropgv}
\max\left(lGV^\delta_{1/n}(f,A,K),rGV^\delta_{1/n}(f,A,K)\right)\\\leq 
\min\left(\overline{GV}^\delta_{1/n}(f,A,K),GV_{1/n}(f,A\cup\{a,b\},K\cup\{a,b\})\right).
\end{multline}
Further, if $0<\delta<\gamma$, then
$\overline{GV}^{\delta}_{1/n}(f,A,K)\leq\overline{GV}^\gamma_{1/n}(f,A,K)$.

We will need some properties of the ``unilateral'' variations.

\begin{lemma}\label{mucontlem}
Let $\emptyset\neq A\subset K\subset[a,b]$ be closed sets, $\{a,b\}\subset K$,
$f:[a,b]\to\R$ continuous and such that $\overline{GV}^\delta_{1/n}(f,A,K)<\infty$
for some $\delta>0$ and $n\geq2$.
Then $v(x):=rGV^\delta_{1/n}(f,A\cap[a,x],K\cap[a,x])$ and
$\tilde v(x):=lGV^\delta_{1/n}(f,A\cap[x,b],K\cap[x,b])$ are continuous
functions on~$K$ such that $v$ is increasing, $\tilde v$ is decreasing,
and
\begin{equation}\label{maxmin}
\max(v(b)-v(a),\tilde v(a)-\tilde v(b))\leq \overline{GV}^\delta_{1/n}(f,A,K).
\end{equation}
\end{lemma}

\begin{proof}
Clearly, $v$ is increasing and $\tilde v$ is decreasing on~$K$.
Also, we easily see that~\eqref{maxmin} holds.
We will only establish the continuity of~$v$ on~$K$, as the case of $\tilde v$ is similar.
To show that $v$ is continuous on $K$, let $x\in K$ be such that $x$
is not a left-hand-side accumulation point of~$A$.  
If $[x-\delta,x)\cap A=\emptyset$,
then $v$ is constant on $[x-\delta,x]$.
If $[x-\delta,x)\cap A\neq\emptyset$, then 
\begin{equation}\label{vvLP} 
v(t)=v(r)+\left(V_{\frac{1}{n-1}}(f,K\cap[r,t])\right)^{\frac{n-1}{n}},
\end{equation}
for all $t\in[r,x]$, where $r:=\max([a,x)\cap A)$ provided $[a,x)\cap A\neq\emptyset$ (the
case when $[a,x)\cap A=\emptyset$ is trivial).
Obviously, we have that $v(t)-v(r)\geq \left(V_{\frac{1}{n-1}}(f,K\cap[r,t])\right)^{\frac{n-1}{n}}$.
The other inequality in~\eqref{vvLP} follows easily from the definition
of the fractional variation together with the fact that $[x-\delta,x)\cap A=\emptyset$.
The unilateral continuity of $v$ at $x$ in this case follows from~\cite[Lemma~3.13]{LP}.
We have shown that $v$ is continuous from the left at all $x\in K$ which
are not left-hand-side accumulation points of $A$.
\par
If $x\in K$ is not right-hand-side accumulation point of $A$, then
either $(x-\delta,x]\cap A=\emptyset$, in which case $v$ is constant on $[x,x+\eta]$ for some
$\eta>0$, or if $(x-\delta,x]\cap A\neq\emptyset$, then~\eqref{vvLP} holds and continuity of $v$
at~$x$ from the right follows again from~\cite[Lemma~3.13]{LP}.
\par
Now suppose that $x\in A$ is a left-hand-side accumulation point of $A$. Fix $\e>0$.
Choose the sequence $([x_i,y_i])^N_{i=1}$ as in~\eqref{Wn} for $rGV^\delta_{1/n}(f,A\cap [a,x],K\cap [a,x])$ 
and $([c^i_j,d^i_j])^{J_i}_{j=1}$ from~\eqref{v1overm}
for $V_{\frac{1}{n-1}}(f,[x_i,y_i]\cap K)$ so that
\begin{equation}\label{ctcond} 
\sum^{N}_{i=1} \bigg( \sum_{j=1}^{J_i} |f(c^i_j)-f(d^i_{j-1})|^{\frac{1}{n-1}}\bigg)^{\frac{n-1}{n}}>v(x)-\e.
\end{equation}
We can assume that $([x_i,y_i])$ and $([c^i_j,d^i_j])$ are ordered
in the natural sense (i.e.\ $y_i<x_{i'}$ whenever $i<i'$; similarly for $c^i_j,d^i_j$).
We can suppose that $d^N_{J_i}<x$ because if
$d^N_{J_i}=x$, then by continuity of $f$ and by the fact that $x$ is a left-hand-side accumulation point of $A$,
we can make $d^N_{J_i}$ slightly smaller (and adjust $y_N$)
without violating~\eqref{ctcond}. Then we have $v(z)>v(x)-\e$ for every ${d^N_{J_i}}<z<x$ and hence $v$
is continuous from the left at $x$.
\par
Now suppose that $x\in A$ is a a right-hand-side point of accumulation of $A$, and fix $\e>0$.
Choose $([x_i,y_i])^N_{i=1}$ and $([c^i_j,d^i_j])^{J_i}_{j=1}$ from~\eqref{Wn} for $rGV_{1/n}^\delta(f,A\cap[x,b])=v(b)-v(x)$
(the equality holds because $x\in A$) so that
\begin{equation}\label{ctcond1} 
\sum^{N}_{i=1} \bigg( \sum_{j=1}^{J_i} |f(c^i_j)-f(d^i_{j-1})|^{\frac{1}{n-1}}\bigg)^{\frac{n-1}{n}}>v(b)-v(x)-\e.
\end{equation}
As before, we can assume that $([x_i,y_i])$ and $([c^i_j,d^i_j])$ are ordered
in the natural sense (see the remark after~\eqref{ctcond}) and that $c^{1}_{1}>x$ (with $x_1>x$). 
Let $y\in(x,x_1)\cap K$. 
Take any $([\tilde x_i,\tilde y_i])^{\tilde N}_{i=1}$ as in~\eqref{Wn} for $rGV^\delta_n(f,A\cap[a,y])=v(y)$.
We have to prove that
\begin{equation}\label{ctcond2} 
\sum^{\tilde N}_{i=1} \left(V_{\frac{1}{n-1}}(f,[\tilde x_i,\tilde y_i]\cap K)\right)^{\frac{n-1}{n}}
\leq v(x)+\e.
\end{equation}
For $i=1,\dots,\tilde N$ take $([\tilde c^i_j,\tilde d^i_j])^{\tilde J_i}_{j=1}$  as
in~\eqref{v1overm} for $V_{\frac{1}{n-1}}(f,[\tilde x_i,\tilde y_i]\cap K)$.
Since $a\leq \tilde x_1<\dots \tilde x_N\leq x_1<\dots<x_N\leq b$, we have 
\begin{equation}\label{ctcond3}
\sum^{\tilde N}_{i=1} \bigg( \sum_{j=1}^{\tilde J_i} |f(\tilde c^i_j)-f(\tilde d^i_{j})|^{\frac{1}{n-1}}\bigg)^{\frac{n-1}{n}}
+\sum^{N}_{i=1} \bigg( \sum_{j=1}^{J_i} |f(c^i_j)-f(d^i_{j})|^{\frac{1}{n-1}}\bigg)^{\frac{n-1}{n}}
\leq v(b).
\end{equation}
The left hand side of~\eqref{ctcond3} is by~\eqref{ctcond1} greater than
\[
\sum^{\tilde N}_{i=1} \bigg( \sum_{j=1}^{\tilde J_i} |f(\tilde c^i_j)-f(\tilde d^i_{j})|^{\frac{1}{n-1}}\bigg)^{\frac{n-1}{n}}
+v(b)-v(x)-\e,\]
and this easily implies~\eqref{ctcond2}.
This concludes the proof of continuity of $v$ on~$K$ since it is clearly an increasing function.

\end{proof}

Now we are ready to define our classes of functions. The following class plays
a similar r\^ole for $n$-times differentiabily as the class $CBV_{1/n}$ from~\cite{LP}
in case of continuous derivatives.

\begin{definition}\label{cbvgdef}
We say that a continuous $f:[a,b]\to\R$ is {\em $CBVG_{1/n}$ for $n\geq2$} provided 
$V_{\frac{1}{n-1}}(f,K_f)<\infty$, and there exist closed $A_m\subset K_f$
$(m\in\mcM\subset\N)$ such that $K_f=\bigcup_{m\in\mcM} A_m$, 
and $GV_{1/n}(f,A_m,K_f)<\infty$ for all $m\in\mcM$.
\end{definition}

It is easy to see that if $f$ is a $CBVG_{1/n}$ function for some $n\geq2$, 
then $f$ has bounded variation.
If $n=2$, then it is not difficult to prove that the class $CBVG_{1/2}$ coincides
with the class $VBG_{\frac{1}{2}}$ from~\cite{D12}.
\par
The analogue of the class $SBV_{1/n}$ from~\cite{LP} in the case of continuous derivatives is
given by the following definition for the case of $n$-times differentiable functions.

\begin{definition}\label{sbvgdef}
We say that a continuous $f:[a,b]\to\R$ is {\em $SBVG_{1/n}$ for $n\geq2$} provided
$V_{\frac{1}{n-1}}(f,K_f)<\infty$ and
there exist closed sets $A_m\subset K_f$, numbers $\delta_m>0$ $(m\in\N)$
such that $K_f=\bigcup_{m\in\N} A_m$ and we have
\begin{enumerate}
\item $\lim_{m\to\infty}\delta_m=0$,
\item $A_m\subset A_{m+1}$ for each $m\in\N$,
\item 
$\lim_{m\to\infty}\overline{GV}^{\delta_m}_{1/n}(f,A_m,K_f)=0$.
\end{enumerate}
\end{definition}

We will need the following auxiliary class:

\begin{definition}\label{sbvgprimedef}
We say that a continuous $f:[a,b]\to\R$ is {\em $\overline{SBVG}_{1/n}$ for $n\geq2$} provided
$V_{\frac{1}{n-1}}(f,K_f)<\infty$ and
there exist closed sets $A_m^k\subset K_f$, numbers $\delta_m^k>0$ $(k,m\in\N)$
such that $K_f=\bigcup_{k,m} A_m^k$ and for each  $k\in\N$ we have
\begin{enumerate}
\item $\lim_{m\to\infty}\delta^k_m=0$,
\item $A^k_m\subset A^k_{m+1}$ for each $m\in\N$,
\item $\lim_{m\to\infty}\overline{GV}^{\delta^k_m}_{1/n}(f,A_m^k,K_f)=0$.
\end{enumerate}
\end{definition}

Since every continuous function on a compact interval is uniformly continuous,
it is easy to see that if~$f$ is $SBVG_{1/n}$ (resp.\ $\overline{SBVG}_{1/n}$ or $CBVG_{1/n}$) 
and~$g$ is Lebesgue equivalent to~$f$, then~$g$ is $SBVG_{1/n}$ (resp. $\overline{SBVG}_{1/n}$ or $CBVG_{1/n}$).
\par
We need the following observation:

\begin{lemma}\label{noprimelem}
Let $f:[a,b]\to\R$. Then $f$ is $SBVG_{1/n}$ if and only if $f$ is $\overline{SBVG}_{1/n}$.
\end{lemma}

\begin{proof} 
Suppose that $f$ is $SBVG_{1/n}$.
Let $A_m$ and $\delta_m$ ($m\in\N$) be as in Definition~\ref{sbvgdef} for~$f$.
We define $A^1_m:=A_m$ and $A^k_m:=\emptyset$ for $k>1$. Similarly,
$\delta^1_m:=\delta_m$ and $\delta^k_m:=1/m$ for $k>1$. Then it is easy
to see that $(A^k_m)$ and $(\delta^k_m)$ satisfy Definition~\ref{sbvgprimedef}
for $f$.
\par
Now, suppose that $f$ is $\overline{SBVG}_{1/n}$,
and let $A^k_m$ and $\delta^k_m$ be as in Definition~\ref{sbvgprimedef}.
By relabeling, we can assume that $A^k_m\neq\emptyset$ for all $m,k\in\N$
(since the case when there exists $k_0\in\N$
such that $\bigcup_m A^k_m=\emptyset$ for all $k\geq k_0$ is simple to handle
using~\eqref{becka}).
By a standard diagonalization argument (using the property~(ii)
from Definition~\ref{sbvgprimedef}), we can also assume that
$\overline{GV}^{\delta^k_m}_{1/n}(f,A^k_m,K_f)\downarrow 0$ when $m\to\infty$
for all $k\in\N$.
Define $A_1:=A^1_1$, $\delta_1:=\delta^1_1$, and $N_1:=1$.
By induction, we will construct closed sets $A_p\subset K_f$, $\delta_p>0$, and $N_p\in\N$.
Suppose that $A_1,\dots,A_{p-1}$ (together with $\delta_i$ and $N_i$ for $i<p$) 
were constructed. Note that for closed $B_1,\dots,B_l\subset K_f$, and $\xi>0$, we have
\begin{equation}\label{becka}
\overline{GV}^{\xi}_{1/n}\bigg(f,\bigcup_{j=1}^l B_j,K_f\bigg)\leq
\sum_{j=1}^l
\overline{GV}^{\xi}_{1/n}(f,B_j,K_f).
\end{equation}
For $i=1,\dots,p$, find $l_i>N_{p-1}$ such that
$\overline{GV}^{\delta^i_{l_i}}(f,A^i_{l_i},K_f)\leq p^{-2}$,
put $A_p:=\bigcup_{i=1}^p A^i_{l_i}$, $\delta_p:=\min(p^{-1},\min_{i=1,\dots,p}\delta^i_{l_i})$,
and using~\eqref{becka} conclude that
\begin{equation}\label{sjednod1}
\overline{GV}^{\delta_p}_{1/n}(f,A_p,K_f)\leq p^{-1}.
\end{equation}
Finally, put $N_p:=1+\max_{i=1,\dots,p}l_i$, and proceed with induction.
Since $K_f=\bigcup_{k,m}A^k_m$, we obtain $K_f=\bigcup_p A_p$.
By construction, it is easy to see that $A_p\subset A_{p+1}$ and $\delta_p\to0$.
Finally,~\eqref{sjednod1} shows that
\[\lim_{p\to\infty}\overline{GV}^{\delta_p}_{1/n}(f,A_p,K_f)=0.\]
\end{proof}

We will also need the following property.

\begin{lemma}\label{ndenselem}
Let $f:[a,b]\to\R$ be a function which is either $CBVG_{1/n}$, $\overline{SBVG}_{1/n}$, or $SBVG_{1/n}$. Then $\lambda(f(K_f))=0$.
\end{lemma}

\begin{proof}
If $f$ is in one of the three classes, then 
$V_{\frac{1}{n-1}}(f,K_f)<\infty$, and thus Lemma~\ref{zerolem}
implies the conclusion. 
\end{proof}

Next lemma contains our basic estimate.

\begin{lemma}\label{nn1odhlem}
Let $n\geq2$, and $f:[a,b]\to\R$ be $(n-1)$-times differentiable with $f^{(n-1)}(\cdot)$ being pointwise
Lipschitz, let $K\subset\{x\in[a,b]:f'(x)=0\}$ be closed with $a,b\in K$, 
and such that $|f(c)-f(d)|=V(f,[c,d])$ for all intervals $(c,d)$ contiguous to $K$
in $[a,b]$, let
\[A\subset\big\{x\in K':|f^{(n-1)}(y)|\leq k|y-x|\text{ for all }y\in B(x,1/m)\big\}\]
be a closed set, where $k,m\in\N$, $a\leq x<x'\leq b$ are such
that $x,x'\in K$, $\{x,x'\}\cap A\neq\emptyset$ and $0<x'-x<1/m$.
Then 
\begin{equation}\label{crodh}
\left( V_{\frac{1}{n-1}}(f,K\cap [x,x'])\right)^{\frac{n-1}{n}}
\leq C_{kn}(x'-x)^{\frac{1}{n}}
\bigg(\sum_{p\in\mcP}(v_p-u_p)\bigg)^{\frac{n-1}{n}},
\end{equation}
where $(u_p,v_p)$ $(p\in\mcP\subset\N)$ are all the intervals contiguous
to~$K$ in~$[x,x']$, and $C_{kn}=k^{\frac{1}{n}}(2n)^{\frac{n-1}{n}}$.
\end{lemma}

\begin{proof}
Without any loss of generality, assume that $x\in A$ (if $x'\in A$, then
work with $f(-\cdot)$ instead).
By Lemma~\ref{Kflem} we have 
\begin{equation}\label{fizero}
f^{(i)}(x)=0\quad\text{ for all }i=1,\dots,n-1.
\end{equation}
Let
\begin{equation}\label{cjdj} 
\begin{split}
&([c_j,d_j])_{j=1}^J\text{ be non-overlapping intervals}\\
&\text{ with }c_j,d_j\in K\cap[x,x']
\text { for all }j\in\{1,\dots,J\}.
\end{split}
\end{equation}
Assume first that $\#(K\cap[x,x'])<2n+1$. 
Then
\begin{equation*}
\begin{split}
|f(c_j)-f(d_j)|&\leq\int^{d_j}_{c_j} |f'(s)|\,ds
\leq \int^{d_j}_{c_j} |f^{(n-1)}(\xi_{n-1})|(x'-x)^{n-2} \,ds\\
&\leq k(x'-x)^{n},
\end{split}
\end{equation*}
where $\xi_i=\xi_i(s)$ is chosen inductively (using~\eqref{fizero}) such that $\xi_1=s$, and
\[|f^{(i-1)}(\xi_{i-1})|=|f^{(i-1)}(\xi_{i-1})-f^{(i-1)}(x)|=|f^{(i)}(\xi_i)|\,|\xi_{i-1}-x|\]
for $i=2,\dots,n-1$.
We obtain $|f(c_j)-f(d_j)|^{\frac{1}{n-1}}\leq k^{\frac{1}{n-1}} (x'-x)^{\frac{n}{n-1}}$,
and it follows that
\begin{equation*}
\begin{split}
\sum_{j=1}^J |f(c_j)-f(d_j)|^{\frac{1}{n-1}}&
\leq k^{\frac{1}{n-1}}
\sum_{j=1}^J (x'-x)^{\frac{n}{n-1}}\\
&\leq k^{\frac{1}{n-1}}(x'-x)^{\frac{1}{n-1}} (2n)\bigg(\sum_{p\in\mcP} (v_p-u_p)\bigg),
\end{split}
\end{equation*}
where we used that $x'-x=\sum_{p\in\mcP}(v_p-u_p)$ and $J\leq 2n$ since $\#(K\cap[x,x'])<2n+1$.
Thus
\[ \bigg(\sum_{j=1}^J |f(c_j)-f(d_j)|^{\frac{1}{n-1}}\bigg)^{\frac{n-1}{n}}
\leq k^{\frac{1}{n}}(2n)^{\frac{n-1}{n}}(x'-x)^{\frac{1}{n}}
\bigg(\sum_{p\in\mcP} (v_p-u_p)\bigg)^{\frac{n-1}{n}},\]
and~\eqref{crodh} holds in this case.
\par
Now assume that $\#(K\cap[x,x'])\geq 2n+1$.
Let $([c_j,d_j])_{j=1}^J$ be as in~\eqref{cjdj}.
Since $\lambda(f(K))=0$ by Lemma~\ref{Sard}, 
Lemma~\ref{basiclem} implies that 
\begin{equation}\label{adding} 
|f(d_j)-f(c_j)|\leq V(f,[c_j,d_j])=\sum_{p\in\mcP_j}V(f,[\gamma_p^j,\delta_p^j])
=\sum_{p\in\mcP_j}|f(\delta_p^j)-f(\gamma_p^j)|,
\end{equation}
where we used the assumptions in the last equality, and where $(\gamma_p^j,\delta_p^j)$ ($p\in\mcP_j\subset\N$) are
all the intervals contiguous to $K\cap [c_j,d_j]$ in $[c_j,d_j]$.
By adding~\eqref{adding} for $j=1,\dots,J$, and using the subadditivity
of $g(t)=t^{\frac{1}{n-1}}$ for $t\geq0$, we obtain
\begin{equation}\label{eee1}
\sum_{j=1}^J|f(c_j)-f(d_j)|^{\frac{1}{n-1}}\leq \sum_{p\in\mcP}|f(v_p)-f(u_p)|^{\frac{1}{n-1}};
\end{equation}
thus the conclusion of the lemma will follow from~\eqref{eee1} once we establish
\begin{equation}\label{withpeq}
\bigg(\sum_{p\in\mcP}|f(v_p)-f(u_p)|^{\frac{1}{n-1}}\bigg)^{\frac{n-1}{n}}
\leq C_{kn}(x'-x)^{\frac{1}{n}}
\bigg(\sum_{p\in\mcP}(v_p-u_p)\bigg)^{\frac{n-1}{n}}.
\end{equation}
\par
In the rest of the proof, we will prove~\eqref{withpeq}.
We need the following definition. 
If $(\alpha,\beta)$ is an interval contiguous to $K\cap[x,x']$ in $[x,x']$, 
then put $r_1(\beta):=\beta$, and
\[ r_{i}(\beta):=\inf\{t\in[r_{i-1}(\beta),x']:\#(K\cap[\beta,t])\geq i\text{ or }t=x'\}\]
for $i=2,\dots,n-2$. 
Similarly, define $l_1(\alpha):=\alpha$, and 
\[ l_{i}(\alpha):=\sup\{t\in[x,l_{i-1}(\alpha)]:\#(K\cap[t,\alpha])\geq i\text{ or }t=x\}\]
for $i=2,\dots,n-2$.
We have the following easy observation:
\begin{itemize}
\item[$(*)$] for each interval 
$(\alpha,\beta)$ contiguous to $K\cap[x,x']$ in $[x,x']$ and each
$i\in\{1,\dots,n-2\}$ there
exists $w_i\in[l_i(\alpha),r_i(\beta)]$ such that $f^{(i)}(w_i)=0$.
\end{itemize}
We will now prove~$(*)$.
If $i=1$, then take $w_1=\alpha$ and~$(*)$ follows. 
Suppose that $1<i<n-1$.
If~$[l_i(\alpha),r_i(\beta)]\cap K'\neq\emptyset$, then~$(*)$ follows by Lemma~\ref{Kflem},
otherwise $\#([l_i(\alpha),r_i(\beta)]\cap K)\geq i+1$ by the choice of $l_i(\alpha),r_i(\beta)$, 
and~$(*)$ follows from Lemma~\ref{Rolle}.
\par
Fix $p\in\mcP$, and suppose that $s\in[u_p,v_p]$. Then
$|f'(s)|=|f'(s)-f'(u_p)|=|f''(\xi_1)|\cdot|s-u_p|$,
and by induction for $k=2,\dots,n-2$, choose $\xi_k\in[l_k(u_p),r_k(v_p)]$ so that
\begin{equation}\label{ee1} 
|f^{(k)}(\xi_{k-1})|=|f^{(k)}(\xi_{k-1})-f^{(k)}(w_k)|
=|f^{(k+1)}(\xi_k)|\cdot|\xi_{k-1}-w_k|,
\end{equation}
where $w_{k}$ is chosen by applying~$(*)$ to $(u_p,v_p)$.
Using~\eqref{ee1}, we obtain
\[|f'(s)|\leq |f^{(n-1)}(\xi_{n-2})|\cdot\prod^{n-2}_{l=1} (r_l(v_p)-l_l(u_p)).\]
Thus we obtain
\begin{equation*}
\begin{split}
|f(u_p)-f(v_p)|&\leq\int^{v_p}_{u_p} |f'(s)|\,ds\leq \int^{v_p}_{u_p} 
|f^{(n-1)}(\xi_{n-2})|\prod^{n-2}_{l=1}(r_l(v_p)-l_l(u_p)) \,ds\\
&\leq k(v_p-x)(r_{n-2}(v_p)-l_{n-2}(u_p))^{n-2}(v_p-u_p).
\end{split}
\end{equation*}
From this we get
$|f(u_p)-f(v_p)|^{\frac{1}{n-1}}\leq (k(v_p-x))^{\frac{1}{n-1}}
(r_{n-2}(v_p)-l_{n-2}(u_p))$
for each $p\in\mcP$, and thus
\begin{equation*}
\bigg(\sum_{p\in\mcP} |f(u_p)-f(v_p)|^{\frac{1}{n-1}}\bigg)^{\frac{n-1}{n}}
\leq k^{\frac{1}{n}}(x'-x)^{\frac{1}{n}}
\bigg(\sum_{p\in\mcP}(r_{n-2}(v_p)-l_{n-2}(u_p))\bigg)^{\frac{n-1}{n}},
\end{equation*}
but since 
$\sum_{p\in\mcP} (r_{n-2}(v_p)-l_{n-2}(u_p))\leq 2n\cdot\sum_{p\in\mcP} (v_p-u_p)$,
where $(u_p,v_p)$ ($p\in\mcP\subset\N$) are all the intervals contiguous to $K$ in $[x,x']$,
we obtain
\[
\bigg(\sum_{p\in\mcP}|f(v_p)-f(u_p)|^{\frac{1}{n-1}}\bigg)^{\frac{n-1}{n}}
\leq
k^{\frac{1}{n}}(2n)^{\frac{n-1}{n}}(x'-x)^{\frac{1}{n}}
\bigg(\sum_{p\in\mcP}(v_p-u_p)\bigg)^{\frac{n-1}{n}}.\]
Thus,~\eqref{withpeq} and also~\eqref{crodh} follow.
\end{proof}

The following lemma contains a sufficient condition guaranteeing that a function 
belongs to the classes $\overline{SBVG}_{1/n}$ and $CBVG_{1/n}$.

\begin{lemma}\label{gwlem}
Suppose that $f:[a,b]\to\R$ is $(n-1)$-times differentiable $($for $n\geq2)$
so that $f^{(n-1)}(\cdot)$ is pointwise-Lipschitz.
Then $f$ is both $\overline{SBVG}_{1/n}$ and $CBVG_{1/n}$.
\end{lemma}

\begin{proof} 
First, note that $V_{\frac{1}{n-1}}(f,K_f)<\infty$ by~\cite[Remark~3.6]{LP}
since $f$ is $C^{n-1}$ by the assumption.
Since $K_f\subset\{x\in[a,b]:f'(x)=0\}$, we have that $f^{(n-1)}(x)=0$ for $x\in K_f'$ by Lemma~\ref{Kflem}.
Denote
\[B^k_m:=\big\{x\in K_f':|f^{(n-1)}(y)|\leq k|y-x|\text{ for all }y\in B(x,1/m)\big\}.\]
It is easy to see that each $B^k_m$ is closed, $\bigcup_{k,m}  B^k_m=K_f'$, and $B^k_m\subset B^k_{m+1}$.
Let $\delta_m^k:=\frac{1}{2m}$.
\par
First, we will show that $\lim_{m\to\infty}\overline{GV}_{1/n}^{\delta^k_m}(f,B^k_m,K_f)=0$.
Fix $k,m\in\N$ with $B^k_m\neq\emptyset$.
Let the collection $([x_i,y_i])_{i=1}^N$ as in~\eqref{Wn} for $A=B_m^k$, $K=K_f$ and $\delta=\delta^k_m$
in the definition of $\overline{GV}^{\delta^k_m}_{1/n}(f,B_m^k,K_f)$.
Let 
\begin{equation}\label{Pis}
\mcP_i:=\{p\in\mcP:(u_p,v_p)\subset[x_i,y_i]\},
\end{equation} 
where
$(u_p,v_p)$ ($p\in\mcP\subset\N$) are all the intervals
contiguous to~$K_f$ in~$[a,b]$.
Applying Lemma~\ref{nn1odhlem} to $[x,x']=[x_{i},y_{i}]$
for a fixed $i=1,\dots,N$,
summing over~$i\in\{1,\dots,N\}$, and using H\"older's inequality (with exponents
$p=n$ and $p'=\frac{n}{n-1}$), we obtain
\begin{equation}\label{Hold}
\begin{split}
\sum^{N}_{i=1}\left(V_{\frac{1}{n-1}}(f,K_f\cap[x_i,y_i])\right)^{\frac{n-1}{n}}
&\leq C_{kn}\sum^{N}_{i=1} (y_i-x_i)^{\frac{1}{n}}\cdot\bigg(\sum_{p\in\mcP_i}(v_p-u_p)\bigg)^{\frac{n-1}{n}}\\
&\leq C_{kn}(b-a)^{\frac{1}{n}}\bigg(\sum^N_{i=1} \sum_{p\in\mcP_i}
(v_p-u_p)\bigg)^{\frac{n-1}{n}}\\
&\leq C_{kn}(b-a)^{\frac{1}{n}}\bigg(\sum_{v_p-u_p\leq\delta^k_m} (v_p-u_p)\bigg)^{\frac{n-1}{n}},
\end{split}
\end{equation}
and since we have $\lim_{m\to\infty}\delta^k_m=0$ for a fixed $k\in\N$,
by~\eqref{Hold} we obtain that 
\begin{equation}\label{hold1}
\lim_{m\to\infty} \overline{GV}_{1/n}^{\delta^k_m}(f,B^k_m,K_f)=0,
\end{equation}
for each $k\in\N$.
\par
Index $K_f\setminus K_f'$ as $\{x_j\}_{j\in\mcJ}$, where $\mcJ\subset\N$
and each $j\in\mcJ$ is even.
Now, define $A^{2k+1}_m:=B^{2k+1}_m$, and take $\tilde\delta^{2k+1}_m:=\delta^{2k+1}_m$.
For $k\in\mcJ$ define 
$A^k_m:=\{x_k\}$ for all $m\in\N$. For each $k\in\mcJ$ find
$\gamma_k>0$ such that $B(x_k,\gamma_k)\setminus\{x_k\}\cap K_f=\emptyset$
and put $\tilde\delta^k_m:=\min(\gamma_k,1/m)$.
For $k\in\N\setminus\mcJ$ which are even, 
put $A^k_m:=\emptyset$ for all $m\in\N$ and $\tilde\delta^k_m:=1/m$.
Using~\eqref{hold1}, it is easy to see that $f$ satisfies Definition~\ref{sbvgprimedef}
with $A^k_m$ and $\tilde\delta^k_m$.
\par
To show that $f$ is $CBVG_{1/n}$,
write each $B^k_m$ as $B^k_m=\bigcup_l B^k_{ml}$, where each $B^k_{ml}$ is closed,
and $\diam(B^k_{ml})<1/m$. 
Fix $k,m,l\in\N$ such that $B^k_{ml}\neq\emptyset$.
Let $([x_i,y_i])^N_{i=1}$, $x_i,y_i\in B^k_{ml}$, be as in~\eqref{Wn} for $f$, $A=B^k_{ml}$,
$K=K_f$ and the definition of $GV_{1/n}(f,B^k_{ml},K_f)$.
If $i\in\{1,\dots,N\}$, then Lemma~\ref{nn1odhlem} applied to $[x,x']=[x_i,y_i]$
shows that
\begin{equation}\label{mezi1}
\begin{split}
\left(V_{\frac{1}{n-1}}(f,[x_i,y_i]\cap K_f)\right)^{\frac{n-1}{n}}
&\leq C_{kn}(y_i-x_i)^{\frac{1}{n}}\bigg(\sum_{i\in\mcP_i}(v_p-u_p)\bigg)^{\frac{n-1}{n}}\\
&\leq C_{kn}(y_i-x_i),
\end{split}
\end{equation}
where $\mcP_i$ is defined as in~\eqref{Pis}.
By summing over $i\in\{1,\dots,N\}$ in~\eqref{mezi1},
we obtain
\begin{equation}\label{perfcase}
\sum_{i=1}^N\left(V_{\frac{1}{n-1}}(f,[x_i,y_i]\cap K_f)\right)^{\frac{n-1}{n}}
\leq C_{kn}\sum_{i=1}^N (y_i-x_i)
\leq C_{kn}(b-a)
<\infty.
\end{equation}
We proved that $GV_{1/n}(f,B^k_{ml},K_f)<\infty$ (for each $k,m,l\in\N$).
If we reorder the sequence $(B^k_{ml})_{k,m,l}$ 
together with the sequence $(\{x\})_{x\in K_f\setminus K_f'}$ 
into a single sequence  which we call $A_m$ (while omitting the empty sets),
where $m\in\mcM\subset\N$, by~\eqref{perfcase} 
we see that $f$ is $CBVG_{1/n}$.
\end{proof}

The following lemma will allow us to construct certain variations,
which play a key r\^ole in establishing differentiability.

\begin{lemma}\label{propvlem}
Let $f:[a,b]\to\R$ be continuous, 
$\emptyset\neq A\subset K\subset[a,b]$ be closed sets, $\{a,b\}\subset K$,
and $\delta>0$.
Suppose that $\overline{GV}^\delta_{1/n}(f,A,K)<\infty$.
Then there exists a continuous increasing function 
$v$ on~$[a,b]$ with $v(a)=0$, $v(b)\leq\overline{GV}_{1/n}^\delta(f,A)$,
and such that for $x\in A$ and $y,z\in K$ with $x\leq y<z<x+\delta$ we have
\begin{equation}\label{odhadv}
|f(y)-f(z)|\leq n^{n-1} (v(z)-v(y))^{n-1}(v(z)-v(x)).
\end{equation}
If $f$ has bounded variation, and $V(f,[\alpha,\beta])=|f(\beta)-f(\alpha)|$
whenever $(\alpha,\beta)$ is an interval contiguous to~$K$, 
then $\lambda(v(K))=0$.
\end{lemma}

\begin{proof} 
Define $v(x):=rGV^\delta_n(f,A\cap[a,x],K\cap[a,x])$ for $x\in K$.
The continuity of~$v$ on~$K$ follows from Lemma~\ref{mucontlem}.
Now, extend $v$ to a continuous function on~$[a,b]$ such that $v$ is continuous and affine
on each $[\alpha,\beta]$, whenever $(\alpha,\beta)$ is an interval contiguous to~$K$ in $[a,b]$.
\par
To prove~\eqref{odhadv}, let $x\leq y<z<x+\delta$ where $x\in A$, and $y,z\in K$. 
By continuity, there is no loss of generality in assuming that $x<y$. 
Fix $\e_0>0$. Choose $([x_i,y_i])^N_{i=1}$, and $([c^i_j,d^i_j])^{J_i}_{j=1}$ such that
\[v(y)=v(x)+\sum^{N}_{i=1} \left( \sum_{j=1}^{J_i} \big|f(d^i_j)-f(c^i_{j-1})\big|^{\frac{1}{n-1}}\right)^{\frac{n-1}{n}}+\e,\]
where $0\leq\e<\e_0$. This can be done because $x\in A$ (see~\eqref{rAdd}). 
We can also assume that $([x_i,y_i])$ and $([c^i_j,d^i_j])$ are
ordered in the natural sense (see the remark after~\eqref{ctcond}).
Then
\begin{equation*}
\begin{split} 
v(z)\geq &v(x)+\sum^{N-1}_{i=1} \bigg( \sum_{j=1}^{J_i} \big|f(d^i_j)-f(c^i_{j-1})\big|^{\frac{1}{n-1}}\bigg)^{\frac{n-1}{n}}\\
&\quad+\bigg( \sum_{j=1}^{J_N} \big|f(d^{N}_j)-f(c^{N}_{j-1})\big|^{\frac{1}{n-1}}+|f(y)-f(z)|^\frac{1}{n-1}\bigg)^{\frac{n-1}{n}}.
\end{split}
\end{equation*}
To simplify the notation, put
$b:=\sum_{j=1}^{J_N} \big|f(d^{N}_j)-f(c^{N}_{j-1})\big|^{\frac{1}{n-1}}$,
and $ a:=b+|f(y)-f(z)|^\frac{1}{n-1}$.
Because of the algebraic identity
\[ u^\frac{n-1}{n}-w^\frac{n-1}{n}=(u-w)\cdot
\frac{\sum^{n-2}_{i=0} u^\frac{i}{n} w^\frac{n-2-i}{n}}{\sum^{n-1}_{i=0} u^\frac{i}{n} w^\frac{n-1-i}{n}}\ ,\]
which is easily seen to be valid for all $u,w\geq0$ with $u+w>0$, we obtain
\[ v(z)-v(y)\geq a^\frac{n-1}{n}-b^\frac{n-1}{n}-\e
=(a-b)\cdot
\frac{\sum^{n-2}_{i=0} a^\frac{i}{n} b^\frac{n-2-i}{n}}{\sum^{n-1}_{i=0} a^\frac{i}{n} b^\frac{n-1-i}{n}}-\e.\]
Because $a\geq b$, we obtain
$\sum^{n-1}_{i=0} a^\frac{i}{n} b^\frac{n-1-i}{n}\leq n a^\frac{n-1}{n}$,
and this together with the inequality $v(z)-v(x)\geq a^\frac{n-1}{n}$ implies
\begin{equation*}
\begin{split}
v(z)-v(y)&\geq |f(y)-f(z)|^\frac{1}{n-1}\cdot\frac{a^\frac{n-2}{n}}{n a^\frac{n-1}{n}}-\e\\
&=\frac{|f(y)-f(z)|^\frac{1}{n-1}}{na^\frac{1}{n}}-\e
\geq\frac{|f(y)-f(z)|^\frac{1}{n-1}}{n(v(z)-v(x))^\frac{1}{n-1}}-\e.
\end{split}
\end{equation*}
To finish the proof of~\eqref{odhadv}, send $\e_0\to0$.
\par
Now, suppose that $f$ has bounded variation and $V(f,[\alpha,\beta])=|f(\beta)-f(\alpha)|$
whenever $(\alpha,\beta)$ is an interval contiguous to $K$. We will show that $\lambda(v(K))=0$.
Let $(c_p,d_p)$ ($p\in\mcP\subset\N$) be all the intervals contiguous to $K$ in $[a,b]$.
First, we will prove that 
\begin{equation}\label{vmprop} 
v(b)-v(a)\leq \sum_{p\in\mcP}
(v(d_p)-v(c_p)).
\end{equation}
To prove~\eqref{vmprop}, fix $\e_0>0$, and
let $([x_i,y_i])^{N}_{i=1}$ be non-overlapping intervals
as in Definition~\ref{gvndef} for $rGV^\delta_{1/n}(f,A,K)$ such
that 
\[v(b)-v(a)=\sum_{i=1}^{N} \left(V_{\frac{1}{n-1}}(f,K\cap[x_i,y_i])\right)^{\frac{n-1}{n}}+\e,\]
for some $0\leq\e<\e_0/3$.
Now for each $i=1,\dots,N$,
find non-overlapping intervals $([c^i_j,d^i_j])_{j=1}^{J_i}$ in $[x_i,y_i]$ 
such that $c^i_j,d^i_j\in K\cap[x_i,y_i]$ and 
\begin{equation}\label{lambda1} 
V_{\frac{1}{n-1}}(f,K\cap[x_i,y_i])
\leq \sum_{j=1}^{J_i} |f(d^i_j)-f(c^i_j)|^{\frac{1}{n-1}}+\bigg(\frac{\e_0}{3N}\bigg)^{\frac{n}{n-1}}.
\end{equation}
For $i\in\{1,\dots,N\}$, we can assume that $d^i_j\leq c^i_{j+1}$ for $j=1,\dots,J^i-1$.
For a moment, fix $i\in\{1,\dots,N\}$. By splitting and regrouping the intervals $([c^i_j,d^i_j])^{J_i}_{j=1}$, we can assume that there is a sequence of finite 
families of intervals $\mcA^i_k$, $k=1,\dots,K_i$, and points $a^i_k\in A\cap[x_i,y_i]$
such that if $(\alpha,\beta)\in\mcA^i_k$, then $\alpha,\beta\in K\cap[x_i,y_i]$, 
if $(\sigma,\tau)\in\mcA^i_l$ where $k<l$, then $\alpha<\beta\leq a^i_l\leq\sigma<\tau$,
\begin{equation}\label{lambda15}
\sum_{j=1}^{J_i}|f(d^i_j)-f(c^i_j)|^{\frac{1}{n-1}}
\leq\sum_{k=1}^{K_i}\sum_{(\alpha,\beta)\in\mcA^i_k}|f(\beta)-f(\alpha)|^{\frac{1}{n-1}},
\end{equation}
and 
\[(c^i_j,d^i_j)
\subset\overline{\bigcup\big\{(\alpha,\beta):(\alpha,\beta)\in\mcA^i_k,k=1,\dots,K_i\big\}}.\]
By Lemma~\ref{basiclem} applied to
$f$ on $[a,b]=[\alpha,\beta]$ for $(\alpha,\beta)\in\mcA^i_k$, and $B=(A\cup\{\alpha,\beta\})\cap[\alpha,\beta]$ (note that
$\lambda(f(A))=0$ since $\lambda(f(K))=0$, and thus $\lambda(f(B))=0$), 
let $(\alpha_l^{\alpha,\beta},\beta_l^{\alpha,\beta})$ ($l\in\{1,\dots,L^i_k\}$) 
be a finite collection of intervals contiguous to $(A\cup\{\alpha,\beta\})\cap[\alpha,\beta]$ in $[\alpha,\beta]$ such that
\begin{equation}\label{lambda2} 
|f(\beta)-f(\alpha)|\leq V(f,[\alpha,\beta])\leq\sum_{l=1}^{L_{ik}}
V\big(f,[\alpha_l^{\alpha,\beta},\beta_l^{\alpha,\beta}]\big)
+
\bigg(\frac{\e_0}{3NK_i|\mcA^i_k|}\bigg)^{{n}}.
\end{equation}
On the set 
$\{ (\alpha_l^{\alpha,\beta},\beta_l^{\alpha,\beta}):
l=1,\dots,L_{ik},k=1,\dots,K_i,(\alpha,\beta)\in\mcA^i_k\}$
define an equivalence
relation $\sim$ in the following way: $I\sim J$ whenever $(\max I,\min J)\cap A=\emptyset$ 
and $(\max J,\min I)\cap A=\emptyset$ (note that one of 
the conditions always holds). By $\mcB^i_q$ ($q=1,\dots,Q_i$) denote the equivalence classes of $\sim$.
We have
\begin{equation}\label{lambda3}
\begin{split}
\sum_{k=1}^{K_i}\sum_{(\alpha,\beta)\in\mcA_k^i}&|f(\beta)-f(\alpha)|^{\frac{1}{n-1}}\\
&\leq \sum_{k=1}^{K_i}\sum_{(\alpha,\beta)\in\mcA_k^i}
\sum_{l=1}^{L_{ik}}\left(V\big(f,[\alpha_l^{\alpha,\beta},\beta_l^{\alpha,\beta}]\big)\right)^{\frac{1}{n-1}}
+\bigg(\frac{\e_0}{3N}\bigg)^{\frac{n}{n-1}}\\
&\leq\sum_{q=1}^{Q_i} \sum_{(\eta,\theta)\in\mcB_q^i}
(V(f,[\eta,\theta]))^{\frac{1}{n-1}}+\bigg(\frac{\e_0}{3N}\bigg)^{\frac{n}{n-1}}.
\end{split}
\end{equation}
By Lemma~\ref{basiclem} and the assumptions, we obtain
\begin{equation}
\label{lambda4}
\sum_{(\eta,\theta)\in\mcB^i_q}
(V(f,[\eta,\theta]))^{\frac{1}{n-1}}
\leq
\sum_{\xi\in\Xi^i_q}|f(\omega_\xi)-f(\Omega_\xi)|^{\frac{1}{n-1}}\leq V_{\frac{1}{n-1}}(f,[\tau^i_q,T^i_q]),
\end{equation}
where $(\omega_\xi,\Omega_\xi)$ (for $\xi\in\Xi^i_q\subset\N$) 
are all the intervals contiguous to $K\cap[\tau^i_q,T^i_q]$
for $\tau^i_q=\inf_{x\in I\in\mcB^i_q} x$, and $T^i_q=\sup_{x\in I\in\mcB^i_q} x$.
By putting the inequalities~\eqref{lambda1},~\eqref{lambda15},~\eqref{lambda2},~\eqref{lambda3} 
and~\eqref{lambda4}, we obtain
\begin{equation}\label{pee1} 
\begin{split}
v(b)-v(a) &\leq \sum_{i=1}^{N}\bigg(
\sum_{j=1}^{J_i} |f(d^i_j)-f(c^i_j)|^{\frac{1}{n-1}}\bigg)^{\frac{n-1}{n}}+\frac{2\e_0}{3}\\
&\leq 
\sum_{i=1}^{N}\bigg(\sum_{k=1}^{K_i}\sum_{(\alpha,\beta)\in\mcA^i_k}
|f(\beta)-f(\alpha)|^{\frac{1}{n-1}}\bigg)^{\frac{n-1}{n}}+\frac{2\e_0}{3}\\
&\leq
\sum_{i=1}^{N}\sum_{q=1}^{Q_i}\bigg(\sum_{(\eta,\theta)\in\mcB^i_q}
(V(f,[\eta,\theta]))^{\frac{1}{n-1}}\bigg)^{\frac{n-1}{n}}+\e_0\\
&\leq
\sum_{i=1}^{N}\sum_{q=1}^{Q_i}
\left(V_{\frac{1}{n-1}}(f,[\tau^i_q,T^i_q]\cap K)\right)^{\frac{n-1}{n}}+\e_0.
\end{split}
\end{equation}
Since $(\tau^i_q,T^i_q)\cap A=\emptyset$, for $x\in[\tau^i_q,T^i_q]\cap K$ we have 
\begin{equation}\label{vugly} 
v(x)=v(z)+\left(V_{\frac{1}{n-1}}(f,K\cap[z,x])\right)^{\frac{n-1}{n}},
\end{equation}
where $z=\max(A\cap[a,\tau^i_q])$,
Lemma~\ref{v1nnul} shows that $\lambda(v([\tau^i_q,T^i_q]\cap K))=0$.
Now, Lemma~\ref{basiclem} applied to $\zeta(x)$
(where $\zeta(x)=\left(V_{\frac{1}{n-1}}(f,K\cap[z,x])\right)^{\frac{n-1}{n}}$
for $x\in[\tau^i_q,T^i_q]\cap K$,
and $\zeta$ is continuous and affine on intervals contiguous
to $[\tau^i_q,T^i_q]\cap K$) implies that 
\[\left(V_\frac{1}{n-1}(f,[\tau^i_q,T^i_q]\cap K)\right)^{\frac{n-1}{n}}
= v(T^i_q)-v(\tau^i_q)\leq\sum_{p\in\mcP^i_q}(v(d_p)-v(c_p)),\]
where $\mcP^i_q=\{p\in\mcP:(c_p,d_p)\subset[\tau^i_q,T^i_q]\}$.
Combining this inequality with~\eqref{pee1},
we get
$v(b)-v(a)\leq \sum_{p\in\mcP}
(v(d_p)-v(c_p))+\e_0$,
and by sending $\e_0\to0$ it follows that~\eqref{vmprop} holds.
Since $v(K)\cap v\big(\bigcup_{p\in\mcP}(c_p,d_p)\big)$ is countable,
we have $\lambda(v(K))=v(b)-v(a)-\lambda\big(\bigcup_{p\in\mcP}(c_p,d_p)\big)=0$.
\end{proof}

By a symmetric argument (this time defining $\tilde v(x):=lGV^\delta_n(f,A\cap[x,b])$), 
we obtain the following:

\begin{lemma}\label{proptildevlem}
Let $f:[a,b]\to\R$ be continuous, 
$\emptyset\neq A\subset K\subset[a,b]$ be closed sets,
$\{a,b\}\subset K$, and $\delta>0$.
Suppose that $\overline{GV}^\delta_{1/n}(f,A,K)<\infty$.
Then 
there exists a continuous decreasing function~$\tilde v$ on~$[a,b]$ with $\tilde v(b)=0$,
$\tilde v(a)\leq \overline{GV}^\delta_{1/n}(f,A,K)$, and such that 
for $x\in A$ and $y,z\in K$ with $x-\delta<z<y\leq x$ we have
\begin{equation}\label{odhadtildev}
|f(y)-f(z)|\leq n^{n-1} (\tilde v(z)-\tilde v(y))^{n-1}(\tilde v(z)-\tilde v(x)).
\end{equation}
If $\{a,b\}\subset K$, $f$ has bounded variation, and $V(f,[\alpha,\beta])=|f(\beta)-f(\alpha)|$
whenever $(\alpha,\beta)$ is an interval contiguous to~$K$, 
then $\lambda(\tilde v(K))=0$.
\end{lemma}

We have the following proposition.

\begin{proposition}\label{gwprop}
Let $f:[a,b]\to\R$ be a ${SBVG}_{1/n}$ function.
Then $f$ is Lebesgue equivalent to an $n$-times differentiable function $\vp$ such that
$\vp^{(i)}(x)=0$ whenever $i\in\{1,\dots,n\}$ and $x\in K_\vp$. Also, $\vp'(x)\neq0$
whenever $x\in[a,b]\setminus K_\vp$.
\par
If $f$ is not constant on any interval, then $\lambda(K_\vp)=0$.
\end{proposition}

\begin{proof}
For a moment, assume that the function $f$ is not constant in any interval.
Lemma~\ref{variacelem} shows that $V(f,[a,b])<\infty$.
Then put $g:=f\circ v_f^{-1}:[0,\ell]\to\R$ (where $\ell:=v_f(b)$). 
Since $\lambda(f(K_f))=0$ by Lemma~\ref{zerolem}, and $v_f(K_f)=K_g$, by
Lemma~\ref{basiclem} we have that 
\[\ell=V(f,[a,b])=\sum_{p\in\mcP} V(g,[u_p,v_p])=\sum_{p\in\mcP}(v_p-u_p),\]
where $(u_p,v_p)$ ($p\in\mcP\subset\N$) are all the intervals contiguous to $K_g$ in $[0,\ell]$,
and thus $\lambda(K_g)=\ell-\lambda\big(\bigcup_{p\in\mcP}(u_p,v_p)\big)=0$.
Putting $G(t)=g\big(\frac{\ell}{b-a}\cdot (t-a)\big)$, $t\in[a,b]$, 
we can assume that $f$ satisfies $\lambda(K_f)=0$ (since $f$ is clearly Lebesgue equivalent to~$G$)
provided $f$ is not constant in any interval.
\par
Let $(A_m)_{m}$ and $(\delta_m)_{m}\subset\R_+\setminus\{0\}$ be the sequences from
Definition~\ref{sbvgdef} for~$f$. Find a monotone sequence $(m_j)_{j\in\N}\subset\N$ such that
$\lim_{j\to\infty}m_j=\infty$, and 
\begin{equation}\label{konveq}
\sum_j j\cdot \overline{GV}^{\delta_{m_j}}_{1/n}(f,A_{m_j},K_f)<\infty.
\end{equation}
Relabel $(A_{m_j})_j$ as $(A_m)_m$, and 
$(\delta_{m_j})_j$ as $(\delta_m)_m$.
Then by~\eqref{konveq} we have
$ \sum_m m\cdot \overline{GV}^{\delta_m}_{1/n}(f,A_m,K_f)<\infty$.
For $x\in[a,b]$ define 
\[v(x):=x+ \sum_m m(v_m(x)-\tilde v_m(x)),\]
where $v_m(x)$ (resp.\ $\tilde v_m(x)$) are the functions~$v$ (resp.\ $\tilde v$)
obtained by applying Lemma~\ref{propvlem} (resp.\ Lemma~\ref{proptildevlem})
to~$f$, $K=K_f$, $A=A_m$, and $\delta=\delta_m$. Note that $v:K_f\to[c,d]$ is a continuous
strictly increasing function, which is onto~$[c,d]$, where $c=v(a)$, $d=v(b)$.
\par
In case that $f$ is not constant on any intervals, by Lemmata~\ref{propvlem}
and~\ref{proptildevlem}, we have that $\lambda(v_m(K_f))=0$, $\lambda(\tilde v_m(K_f))=0$
for each $m\in\N$. Also, $\lambda(K_f)=0$, and thus Lemma~\ref{propnlem} shows
that 
\begin{equation}\label{lambda0}
\lambda(v(K_f))=0.
\end{equation}
\par
For $x\in K_f$, we will show that for each $\e>0$ there exists $\delta>0$ 
such that if $x\leq y<z<x+\delta$ or $x-\delta<z<y\leq x$, and $y,z\in K_f$,
then
\begin{equation}\label{varat1}
|f(y)-f(z)|\leq \e|v(z)-v(y)|^{n-1}|v(z)-v(x)|.
\end{equation}
To prove~\eqref{varat1}, fix $x\in K_f$, and $\e>0$. 
Find $m_0\in\N$ such that $x\in A_{m_0}$ (and thus $x\in A_m$ for all $m\geq m_0$),
and pick $m>m_0$ such that $\frac{n^{n-1}}{m^n}<\e$. Define $\delta:=\delta_m$.
Let $y,z\in K_f$ be such that $x<y<z<x+\delta$.
Then~\eqref{odhadv} implies that
$|f(y)-f(z)|\leq n^{n-1} (v_m(z)-v_m(y))^{n-1}(v_m(z)-v_m(x))$.
But since $m(v_m(\tau)-v_m(\sigma))\leq v(\tau)-v(\sigma)$ for all $a\leq\sigma<\tau\leq b$,
we obtain
\[ m^n |f(y)-f(z)|\leq n^{n-1} (v(z)-v(y))^{n-1}(v(z)-v(x)).\]
By the choice of $m$ we have
$|f(y)-f(z)|\leq \e(v(z)-v(y))^{n-1}(v(z)-v(x))$.
By continuity, the above argument shows that~\eqref{varat1} holds also for $y,z\in K_f$
such that $x=y<z<x+\delta$.
Finally, by using~\eqref{odhadtildev} (instead of~\eqref{odhadv}) in the previous argument,
we obtain~\eqref{varat1} for $x-\delta<z<y\leq x$
with $y,z\in K_f$.
\par
We will define  $F:[c,d]\to\R$ as
\begin{equation}\label{defF}
F(x):=\begin{cases}
f\circ v^{-1}(x)&\text{ for }x\in v(K_f),\\
H_{\alpha,\beta}(x)&\text{ for }x\in(\alpha,\beta	),
\end{cases}
\end{equation}
whenever $(\alpha,\beta)$ is an interval contiguous to $v(K_f)$ in $[c,d]$,
and $H=H_{\alpha,\beta}$ is is chosen by applying Lemma~\ref{homeolem}
to $\alpha,\beta, A=f\circ v^{-1}(\alpha), B=f\circ v^{-1}(\beta)$.
It follows that $F$ is Lebesgue equivalent to~$f$,
and $F$ is $n$-times differentiable at all $x\in[c,d]\setminus v(K_f)$
(by Lemma~\ref{homeolem}).
To prove that $F$ is $n$-times differentiable, 
it remains to show that $F^{(i)}(x)=0$ for all $x\in v(K_f)$, $i=1,\dots,n$.
\par
Now,~\eqref{varat1}
implies that for each $x\in v(K_f)$ and for each $\e>0$ there exists $\delta>0$ such that
\begin{equation}\label{eqa1}
|F(y)-F(z)|\leq \e|y-z|^{n-1}|x-z|,
\end{equation}
whenever $x\leq y<z<x+\delta$ or $x-\delta<z<y\leq x$, and $y,z\in v(K_f)$. 
Fix $x\in v(K_f)$.
First, we will show that for each $x\in v(K_f)$ we have that
\begin{equation}\label{pteq}
\begin{split}
\text{ for each }\e>0\text{ there exists }\delta&>0\text{ such that the inequality}\\
|F^{(n-1)}(t)|\leq\e|t-x|\text{ holds for all }&t\in[c,d]\setminus v(K_f)\text{ with }|x-t|<\delta.
\end{split}
\end{equation}
To prove~\eqref{pteq}, if $x$ is not a right-hand-side accumulation 
point of the set $v(K_f)$, then~\eqref{pteq} follows from~\eqref{defF} and the fact
that $H_{\alpha,\beta}^{(n)}(x)=0$ (where $x=\alpha$, and $(\alpha,\beta)$ is
the corresponding interval contiguous to~$v(K_f)$).
If $x$ is a right-hand-side accumulation point of $v(K_f)$,
let $\e>0$ be given, and choose $\delta>0$ such that~\eqref{eqa1} holds
and with $x+\delta\in v(K_f)$.
Let $(\alpha,\beta)$ be an interval contiguous to $v(K_f)$ in $[c,d]$ with $(\alpha,\beta)\subset[x,x+\delta]$, and
let $t\in(\alpha,\beta)$ (the case $t<x$ being treated symmetrically).
Let $l_1(x):=\xi(x-\alpha)$ for $x\in\big[\alpha,\frac{\alpha+\beta}{2}\big]$,
and $l_2(x):=\xi(\beta-x)$ for $x\in\big[\frac{\alpha+\beta}{2},\beta\big]$,
where $\xi=C\cdot\frac{|f\circ v^{-1}(\beta)-f\circ v^{-1}(\alpha)|}{(\beta-\alpha)^n}$,
and $C=C_{\alpha,\beta}$ comes from condition~(iii) of Lemma~\ref{homeolem} applied on~$[\alpha,\beta]$.
By~\eqref{eqa1}, we have that 
\begin{equation*}
\begin{split}
l_i\bigg(\frac{\alpha+\beta}{2}\bigg)
&=C\cdot\frac{|f\circ v^{-1}(\beta)-f\circ v^{-1}(\alpha)|}{(\beta-\alpha)^n}\cdot
\frac{\beta-\alpha}{2}\\
&=\frac{C}{2}\cdot\frac{|f\circ v^{-1}(\beta)-f\circ v^{-1}(\alpha)|}{(\beta-\alpha)^{n-1}}
\leq C\cdot\e\cdot \bigg(\frac{\alpha+\beta}{2}-\alpha\bigg),
\end{split}
\end{equation*}
for $i=1,2$, and this inequality together with the equalities $l_1(\alpha)=l_2(\beta)=0$,
and condition~(iii) from Lemma~\ref{homeolem}
easily implies that 
$|F^{(n-1)}(t)|=|H_{\alpha,\beta}^{(n-1)}(t)|\leq\min(l_1(t),l_2(t))\leq C\e(t-x)$
for all $t\in(\alpha,\beta)$. To see this, we use the fact that if two continuous affine functions $a_1,a_2:[a,b]\to\R$
satisfy $a_1(a)\leq a_2(a)$, and $a_1(b)\leq a_2(b)$, then $a_1(t)\leq a_2(t)$ for all $t\in[a,b]$.
We apply this fact to $a_1(t)=l_1(t)$ for $t\in\big[\alpha,\frac{\alpha+\beta}{2}\big]$ (resp.\ $a_1(t)=l_2(t)$
for $t\in\big[\frac{\alpha+\beta}{2},\beta\big]$), and $a_2(x)=C\e(t-x)$.
Similarly for $(\alpha,\beta)\subset[x-\delta,x]$.
\par
Let $\e>0$, and let $\delta>0$ be as in~\eqref{pteq}. 
It follows easily by induction (using~\eqref{pteq})
that 
\begin{equation}\label{derom}
|F^{(i)}(t)|\leq C\e|t-x|\quad\text{ for all }t\in(x,x+\delta)\setminus v(K_f)\text{ and }i=1,\dots,n-1.
\end{equation}
To show~\eqref{derom}, let $t':=\max (v(K_f)\cap[x,t])$.
Then (using the Mean Value Theorem and the fact
that on $[t',t]$ it holds that $F=H_{\alpha,\beta}$ for some $\alpha,\beta$) we obtain 
$|F^{(i)}(t)|=|F^{(i)}(t)-F_+^{(i)}(t')|\leq|F^{(i+1)}(\xi_{i+1})|\cdot|t-t'|
\leq \dots\leq|F^{(n-1)}(\xi_{n-1})|\cdot|t-t'|^{n-1-i}$,
where $\xi_{j}\in[t',t]$, and~\eqref{derom} easily follows.

Using~\eqref{eqa1},~\eqref{derom}, and the fact that $F^{(i)}_+(\alpha)=H^{(i)}_{\alpha,\beta}(\alpha)=0$ for $i=1,\dots,n$,
by a simple induction argument we obtain that 
\begin{equation}\label{dernul}
F^{(i)}(x)=0\quad\text{for all } x\in v(K_f),\text{ and } i=1,\dots,n-1.
\end{equation}
To prove this, note that the case $i=1$ follows directly from~\eqref{eqa1}. 
For $i=2$, $F'(t)-F'(x)=0$ provided $t\in v(K_F)$ and given $\e>0$ then
for $t\in(x-\delta,x+\delta)\setminus v(K_f)$ (where $\delta$ is chosen so that~\eqref{eqa1} and~\eqref{derom} hold) we have
$|F'(t)-F'(x)|=|F'(t)|\leq C\e |t-x|$, and thus $F''(x)=0$. Similarly for higher $i$'s.
\par
To finish the proof of differentiability of $F$, we will show that $F^{(n)}(x)=0$ for each $x\in v(K_f)$.
But since $F^{(n-1)}(w)=0$ for all $w\in v(K_f)$, \eqref{pteq} together with~\eqref{dernul} easily imply this assertion.
\par
If $f$ is not constant on any interval, then~\eqref{lambda0} implies that
$\lambda(K_F)=\lambda(v(K_f))=0$. 
\par
By~\eqref{defF} and by the property~(i)
of Lemma~\ref{homeolem}, it follows that $F'(x)\neq0$ for all $x\in[c,d]\setminus K_F$.
It follows that there exists a linear homeomorphism $\eta:[a,b]\to[c,d]$, which is onto.
Now, it clearly suffices to put $\vp:=F\circ\eta$.
\end{proof}

\begin{proposition}\label{ptwprop}
Let $f:[a,b]\to\R$ be a $CBVG_{1/n}$ function.
Then $f$ is Lebesgue equivalent to an $(n-1)$-times differentiable function $\phi$ 
such that $\phi^{(n-1)}(\cdot)$ is pointwise Lipschitz, $\phi^{(i)}(x)=0$ for
all $x\in K_\phi$, $i=1,\dots,n-1$, and $\phi'(x)\neq0$ whenever
$x\in[a,b]\setminus K_\phi$.
\par
If $f$ is not constant on any interval, then $\lambda(K_\phi)=0$.
\end{proposition}

\begin{proof} 
It is similar to the proof of Proposition~\ref{gwprop}, so we will only sketch it.
If $f$ is not constant on any interval, similarly as in the proof of Proposition~\ref{gwprop}
we can assume that $\lambda(K_f)=0$.
\par
Let $(A_m)_{m\in\mcM}$ be the sets from Definition~\ref{cbvgdef} for $f$.
If $\mcM$ is finite, then use $\tilde A_m=K_f$ and $\tilde\mcM=\N$ instead
of $A_m$ and $\mcM$ (it is easy to see that in this case $GV_{1/n}(f,K_f,K_f)<\infty$).
Find a sequence of $a_m>0$ such that 
\[\sum_{m\in\mcM} a_m\cdot GV_{1/n}(f,A_{m}\cup\{a,b\},K_f)<\infty.\]
For $x\in[a,b]\cap K_f$ define $v(x)$ as 
$x+\sum_{m\in\mcM} a_m\cdot (v_m(x)-\tilde v_m(x))$,
where $v_m$ (resp.\ $\tilde v_m$) are the functions
obtained from Lemma~\ref{propvlem} (resp.\ Lemma~\ref{proptildevlem}) applied
to $f$ and $\delta=b-a$ (since it is easy to see that
$\overline{GV}^{(b-a)}_{1/n}(f,A_m,K_f)\leq 2\,GV_{1/n}(f,A_m\cup\{a,b\},K_f)$;
see also~\eqref{basicpropgv}).
As in the proof of Proposition~\ref{gwprop}, 
we obtain that $v$ is a continuous strictly increasing function.
We have that for each $x\in K_f$ there exists $m\in\mcM$ such that
\begin{equation}\label{od123} 
|f(y)-f(z)|\leq C_m |v(z)-v(y)|^{n-1}|v(z)-v(x)|
\end{equation}
for all $y,z\in K_f$ with $x\leq y<z$ or $z<y\leq x$. 
This follows from~\eqref{odhadv}, and~\eqref{odhadtildev} in a similar way as~\eqref{varat1}
in the proof of Proposition~\ref{gwprop}. 
\par
Define a new function $F:[a,b]\to\R$ as
\begin{equation}\label{defF1}
F(x):=\begin{cases}
f\circ v^{-1}(x)&\text{ for }x\in v(K_f),\\
H_{\alpha,\beta}(x)&\text{ for }x\in(\alpha,\beta),
\end{cases}
\end{equation}
whenever $(\alpha,\beta)$ is an interval contiguous to $v(K_f)$ in $[c,d]$,
and $H=H_{\alpha,\beta}$ is is chosen by applying Lemma~\ref{homeolem}
to $\alpha,\beta, A=f\circ v^{-1}(\alpha), B=f\circ v^{-1}(\beta)$.
It follows that $F$ is $n$-times differentiable at all $x\in[c,d]\setminus v(K_f)$.
It remains to show that $F^{(i)}(x)=0$ for all $x\in v(K_f)$, $i=1,\dots,n-1$,
and that $F^{(n-1)}$ is pointwise Lipschitz at points $x\in v(K_f)$.
\par
From~\eqref{od123}, we have that for each $x\in v(K_f)$ there exists $C_x>0$ such that
\begin{equation}\label{pteqa}
|F(y)-F(z)|\leq C_x |z-y|^{n-1}|z-x|
\end{equation}
for $y,z\in v(K_f)$ with $z<y\leq x$ or $x\leq y<z$. 
From this we obtain that
\begin{equation}\label{ptw} 
|F^{(n-1)}(t)|\leq C_x'|t-x| 
\end{equation}
for all $t\in [c,d]\setminus v(K_f)$. By induction, we obtain that $F^{(i)}(x)=0$ for
all $i=1,\dots,n-1$ (since~\eqref{pteqa} together with~\eqref{ptw} imply that
$F'(x)=0$, and then~\eqref{ptw} easily implies that $F^{(i)}(x)=0$ for $i=2,\dots,n-1$). 
Now~\eqref{ptw} (together with the fact that $F^{(n-1)}(x)=0$ for each
$x\in v(K_f)$) easily implies that $F^{(n-1)}$ is pointwise Lipschitz at all points of $v(K_f)$.
\par
If $f$ is not constant on any interval, then as in Proposition~\ref{gwprop},
we establish that $\lambda(K_F)=\lambda(v(K_f))=0$. 
\par
By~\eqref{defF1} and by the property~(i)
of Lemma~\ref{homeolem}, it follows that $F'(x)\neq0$ for all $x\in[c,d]\setminus K_F$.
It follows that there exists a linear homeomorphism $\eta:[a,b]\to[c,d]$, which is onto.
Now, it clearly suffices to put $\phi:=F\circ\eta$.
\end{proof}

\section{Main result}\label{mainsec}

The case $n=2$ is handled in~\cite{D12}. The following main theorem gives a slightly
different characterization in that case (see Introduction).

\begin{theorem}\label{mainthm}
Let $f:[a,b]\to\R$ be continuous, $n\geq2$, and $n\in\N$.
Then the following are equivalent:
\begin{enumerate}
\item $f$ is Lebesgue equivalent to a function $g$ which is $n$-times differentiable.
\item  $f$ is Lebesgue equivalent to a function $g$ which is $n$-times differentiable
and such that $g^{(i)}(x)=0$ whenever $i\in\{1,\dots,n\}$ and $x\in K_g$,
and $g'(x)\neq0$ whenever $x\in [a,b]\setminus K_g$.
\item $f$ is Lebesgue equivalent to a function $g$ which is $(n-1)$-times differentiable
and such that $g^{(n-1)}(\cdot)$ is pointwise-Lipschitz.
\item $f$ is $\overline{SBVG}_{1/n}$.
\item $f$ is $CBVG_{1/n}$.
\item $f$ is $SBVG_{1/n}$.
\end{enumerate}
\end{theorem}

\begin{proof}
The implications (ii)$\implies$(i), and (i)$\implies$(iii) are trivial.
The implications (iii)$\implies$(iv), and (iii)$\implies$(v) follow from Lemma~\ref{gwlem}. 
The implication (vi)$\implies$(ii) follows from Proposition~\ref{gwprop},
and the implication (v)$\implies$(iii) from Proposition~\ref{ptwprop}.
Finally, the implication (i)$\implies$(vi) follows from Lemmata~\ref{gwlem} and~\ref{noprimelem};
and the implication (iv)$\implies$(vi) from Lemma~\ref{noprimelem}.
\end{proof}

We have the following corollary:

\begin{corollary}
Let $f:[a,b]\to\R$, $n\geq 2$, $n\in\N$.
Then $f$ is $CBVG_{1/n}$ if and only if $f$ is $SBVG_{1/n}$.
\end{corollary}

The following example shows that for each $n\geq2$ there exists a continuous
function $f:[0,1]\to\R$
such that $f$ is $CBVG_{1/n}$ (and thus $f$ is Lebesgue equivalent to an $n$-times
differentiable function by Theorem~\ref{mainthm}), but $V_{1/n}(f,K_f)=\infty$
(and thus $f$ is not Lebesgue equivalent to any $C^n$ function by
the results of~\cite{LP}). It is a simplified version of~\cite[Example~8.3]{LP}.

\begin{example}\label{difex}
Let $n\geq2$ be an integer. Let $a_m\subset(0,1)$ be such that $a_m\downarrow0$. Define
$f(a_{2m})=m^{-n}$ and $f(_{2m-1})=0$ for all $m=1,2,\dots$, $f(0)=f(1)=0$,
and extend $f$ onto $[0,1]$ such that it is continuous and affine
on each interval contiguous to $K=\{0,1\}\cup\{a_m:m\in\N\}$.
Then $K_f=K$, $f$ is $CBVG_{1/n}$, but $V_{1/n}(f,K_f)=\infty$.
\end{example}

\begin{proof}
Obviously, the function $f$ is continuous, has bounded variation,
and $V_{\frac{1}{n-1}}(f,K_f)<\infty$. Also, it is
easy to see that $f$ is $CBVG_{1/n}$ (using $A_1=\{0,1\}$
and $A_m=\{a_{2m-2},a_{2m-3}\}$ for $m=2,3,\dots$). On the other hand,
\[ V_{1/n}(f,K_f)\geq \sum_{m\in\N} |f(a_{2m})-f(a_{2m-1})|^{\frac{1}{n}}
=\sum_m \frac{1}{m}=+\infty.\]
\end{proof}

The following theorem characterizes the situation when we require
the first derivative to be nonzero almost everywhere.

\begin{theorem}\label{mainnonzerothm}
Let $f:[a,b]\to\R$ be continuous, $n\geq2$, and $n\in\N$.
Then the following are equivalent:
\begin{enumerate}
\item $f$ is Lebesgue equivalent to a function $g$ which is $n$-times differentiable
and such that $g'(x)\neq0$ for a.e.\ $x\in[a,b]$.
\item  $f$ is Lebesgue equivalent to a function $g$ which is $n$-times differentiable,
such that $g^{(i)}(x)=0$ whenever $i\in\{1,\dots,n\}$, $x\in K_g$,
$g'(x)\neq0$ whenever $x\in [a,b]\setminus K_g$, and such that $\lambda(K_g)=0$.
\item $f$ is Lebesgue equivalent to a function $g$ which is $(n-1)$-times differentiable
and such that $g^{(n-1)}(\cdot)$ is pointwise-Lipschitz,
$g'(x)\neq0$ for a.e.\ $x\in[a,b]$.
\item $f$ is $\overline{SBVG}_{1/n}$ and $f$ is not constant on any interval.
\item $f$ is $CBVG_{1/n}$ and $f$ is not constant on any interval.
\item $f$ is $SBVG_{1/n}$ and $f$ is not constant on any interval.
\end{enumerate}
\end{theorem}

\begin{proof}
The proof is similar to the proof of Theorem~\ref{mainthm} while
we also use the fact that $g'(x)\neq0$ for a.e.\ $x\in[a,b]$ implies that
$g$ is not constant on any interval, the fact that being nonconstant
on any interval is invariant with respect to the Lebesgue equivalence,
and the corresponding assertions in Propositions~\ref{gwprop} and~\ref{ptwprop}.
\end{proof}

\section{Generalized Zahorski Lemma}\label{Zahorski}

Our methods yield the following theorem, which can be viewed
as a generalization of Zahorski's Lemma; see e.g.~\cite{Zah1}
or~\cite[p.~27]{GNW}.

\begin{theorem}\label{zahthm}
Let $K\subset [a,b]$ be a closed set, $n\geq2$, $n\in\N$. 
Then the following are equivalent:
\begin{enumerate}
\item There exists an $n$-times
differentiable homeomorphism $h$ of $[a,b]$ onto itself such
that 
\begin{equation*}\label{imeq}
K=h(\{x\in[a,b]:h^{(i)}(x)=0\ \forall i=1,\dots,n\}).
\end{equation*}
\item 
$V_{\frac{1}{n-1}}(id,K)<\infty$ and
there exist closed sets $A_m\subset[a,b]$ $(m\in\mcM\subset\N)$
such that $K=\bigcup_m A_m$ and $GV_{1/n}(id,A_m,K)<\infty$ for all $m\in\mcM$.
\item 
$V_{\frac{1}{n-1}}(id,K)<\infty$ and
there exist closed sets $A_m\subset[a,b]$ and numbers $\delta_m>0$ $(m\in\N)$ such that
$K=\bigcup_m A_m$, $A_m\subset A_{m+1}$ for all $m\in\N$, and $\lim_{m\to\infty}\overline{GV}_{1/n}^{\delta_m}(id,A_m,K)=0$.
\end{enumerate}
\end{theorem}

\begin{proof} 
Since the proof is similar to the considerations above, we will
only sketch it. Without any loss of generality, we can assume that $\{a,b\}\subset K$.
The proof of the implication (i)$\implies$(ii) 
is similar to the proof that every $(n-1)$-times
differentiable function $f$ such that $f^{(n-1)}$ is pointwise Lipschitz, 
is $CBVG_{1/n}$ (see Lemma~\ref{gwlem}).
\par
To prove that (ii)$\implies$(iii), note that an argument similar to 
the proof of Proposition~\ref{ptwprop} shows that there exists a homeomorphism
$\vp$ of $[a,b]$ onto itself which is $(n-1)$-times differentiable
with $\vp^{(n-1)}$ being pointwise Lipschitz, and such that
\[K=\vp(\{x\in[a,b]:\vp^{(i)}(x)=0\text{ for all }i=1,\dots,n-1\}).\]
By the argument of the proof of Lemma~\ref{gwlem},
we find a decomposition $\tilde A^k_m$ with
some $\tilde\delta^k_m>0$ for $\vp^{-1}(K)$ such that
$\tilde A^k_m\subset\tilde A^k_{m+1}$ for all $m,k\in\N$,
\[\lim_{m\to\infty}\overline{GV}^{\tilde\delta_m^k}_{1/n}(id,\tilde A_m^k,\vp^{-1}(K))=0\] for
each $k$, and then
using the diagonal argument from Lemma~\ref{noprimelem},
we find $\tilde A_m$ and $\tilde \delta_m$ such that $\vp^{-1}(K)=\bigcup_m\tilde A_m$,
and 
\[\lim_{m\to\infty}\overline{GV}^{\tilde\delta_m}_{1/n}(id,\tilde A_m,\vp^{-1}(K))=0.\]
Now, we put $A_m:=\vp(\tilde A_m)$, and find suitable $\delta_m>0$ using the uniform continuity
of $\vp^{-1}$. This shows that~(iii) holds. 
\par
Finally, to show that (iii)$\implies$(i), one can use a similar construction
as in the proof of Proposition~\ref{gwprop}.
\end{proof} 

\begin{remark}
It is also not very difficult to see that
in the previous theorem, we can replace~$(i)$ with
\begin{itemize}
\item[(i')] 
There exists an $n$-times
differentiable homeomorphism $h$ of $[a,b]$ onto itself such
that $h^{-1}$ is absolutely continuous, $h^{-1}(K)=0$, and 
\begin{equation*}\label{imeq1}
K=h(\{x\in[a,b]:h^{(i)}(x)=0\ \forall i=1,\dots,n\}).
\end{equation*}
\end{itemize}
The proof uses ideas from the proof of Theorem~\ref{mainnonzerothm}.
Let us only remark that it is well known that
$h^{-1}$ can be taken absolutely continuous in the classical Zahorski's Lemma;
see e.g{.}\ \cite{Br}.
\end{remark}

\section*{Acknowledgment} 
The author would like to thank Professor Lud\v{e}k Zaj\'\i\v{c}ek
for several remarks that led to considerable improvements of the presentation.

\end{document}